\documentclass[12pt]{amsart}
\usepackage{latexsym, amssymb, amsmath, amscd}
 \usepackage{eucal}

    \newtheorem{rema}{Remark}[section]
    \newtheorem{propo}[rema]{Proposition}

   \newtheorem{theo}[rema]{Theorem}
   \newtheorem{def-theo}[rema]{Definition-Theorem}
 
   \newtheorem{defi}[rema]{Definition}
    \newtheorem{lemma}[rema]{Lemma}
    \newtheorem{corol}[rema]{Corollary}
     
  \newtheorem{rmk}[rema]{Remark}

	\newcommand{\nno}{\nonumber}

	\newcommand{\p}{\partial}

 \newcommand{\pf}{{\it Proof:}\hspace{2ex}}
 
 \newcommand{\epfv}{\hspace{1em}$\Box$\vspace{1em}}

\newcommand{\bQ}{{\mathbb Q}}
\newcommand{\bN}{{\mathbb N}}
\newcommand{\bT}{{\mathbb T}}
\newcommand{\bbT}{{\bar{\mathbb T}}}

\newcommand{\cT}{{\mathcal T}}

\newcommand{\cS}{{\mathcal S}}

\newcommand{\cD}{{\mathcal D}}
\newcommand{\cH}{{\mathcal H}}

\newcommand{\cNcs}{{${\mathcal N}$CS} }
\newcommand{\cNsf}{{{\mathcal N}Sym}}
\newcommand{\cQf}{{{\mathcal Q}Sym}}
\newcommand{\Oft}{ \Omega_{F_t}  }
\newcommand{\OTW}{{\Omega_\bT^W} }
\newcommand{\cSft}{ {\mathcal S}_{F_t} }
\newcommand{\cAft}{ {\mathcal A}_{F_t} }

\newcommand{\cDz}{{\mathcal D \langle  z \rangle}}
\newcommand{\cDaz}{{{\mathcal D}^{[\alpha]} \langle  z \rangle}}

\newcommand{\cDzz}{{\mathcal D \langle \langle z \rangle\rangle}}
\newcommand{\cDkzz}{{\mathcal D_K \langle \langle z \rangle\rangle}}
\newcommand{\cDtzz}{{\mathcal D_t \langle \langle z \rangle\rangle}}
\newcommand{\cDrzz}{{\mathcal D er\langle \langle z \rangle\rangle}}
\newcommand{\cDrtzz}{{\mathcal D er_t \langle \langle z \rangle\rangle}}
\newcommand{\cDrkzz}{{\mathcal D er_K \langle \langle z \rangle\rangle}}

\newcommand{\cDrazz}{{\cD er^{[\alpha]}\langle \langle z \rangle \rangle}} 
 
\newcommand{\cDazz}{{\cD^{[\alpha]}\langle \langle z \rangle \rangle}}

\newcommand{\cDrtazz}{{\cD er^{[\alpha]}_t\langle \langle z \rangle \rangle}} 
 
\newcommand{\cDtazz}{{\cD^{[\alpha]}_t\langle \langle z \rangle \rangle}}

\newcommand{\ataz}{{\mathbb A_t^{[\alpha]}\langle \langle z\rangle\rangle}}
\newcommand{\gtaz}{{\mathbb G_t^{[\alpha]}\langle \langle z\rangle\rangle}}

\newcommand{\klam}{{K\langle \Lambda \rangle}}

\newcommand{\kz}{{K\langle z \rangle}}
\newcommand{\kzz}{{K\langle \langle z \rangle\rangle}}
\newcommand{\ktz}{{K[t]\langle z \rangle}}

\newcommand{\kttzz}{{K[[t]]\langle \langle z \rangle\rangle}}

\newcommand{\BQ}{\begin{eqnarray}}
\newcommand{\EQ}{\end{eqnarray}}
\newcommand{\BQn}{\begin{eqnarray*}}
\newcommand{\EQn}{\end{eqnarray*}}

\newcommand{\lb}{\left[}
\newcommand{\rb}{\right]}
\newcommand{\lp}{\left(}
\newcommand{\rp}{\right)}

\newcommand{\fr}{\frac}
\newcommand{\pz}{\frac{\p}{\p z}}
\newcommand{\wtilde}{\widetilde}

\title[\cNcs Systems over Associative Algebras]
{Noncommutative Symmetric Systems over Associative Algebras}

    \author{Wenhua Zhao}      
    \date{\today}
    \begin{document}

\begin{abstract}
This paper is the first of a sequence papers 
(\cite{GTS-II}--\cite{GTS-V})
on the {\it \cNcs $(\text{noncommutative symmetric})$ systems}
over differential operator algebras 
in commutative or noncommutative variables 
(\cite{GTS-II}); the \cNcs systems over
the Grossman-Larson Hopf algebras 
(\cite{GL},\,\cite{F})
of labeled rooted trees (\cite{GTS-IV});
as well as their connections and 
applications to the inversion problem 
(\cite{BCW},\,\cite{E}) 
and specializations of NCSFs 
(\cite{GTS-III},\,\cite{GTS-V}).
In this paper, inspired by 
the seminal work 
\cite{G-T} on NCSFs
(noncommutative symmetric functions), 
we first formulate the notion 
{\it \cNcs systems} 
over associative $\bQ$-algebras.  
We then prove some results for 
\cNcs systems in general; 
the \cNcs systems over 
bialgebras or Hopf algebras; 
and the universal \cNcs system 
formed by the generating functions 
of certain NCSFs in \cite{G-T}. 
Finally, we review some of the main results 
that will be proved in the 
followed papers \cite{GTS-II}, 
\cite{GTS-IV} and \cite{GTS-V} 
as some supporting examples 
for the general discussions 
given in this paper.
\end{abstract}

   

 \bibliographystyle{alpha}
    \maketitle


\renewcommand{\theequation}{\thesection.\arabic{equation}}
\renewcommand{\therema}{\thesection.\arabic{rema}}
\setcounter{equation}{0}
\setcounter{rema}{0}
\setcounter{section}{0}

\section{\bf Introduction}\label{S1}

Let $K$ be any unital commutative 
$\bQ$-algebra and $A$ a 
unital associative but not necessarily 
commutative $K$-algebra. 
Let $t$ be a formal central parameter, 
i.e. it commutes with all elements of $A$, 
and $A[[t]]$ the $K$-algebra 
of formal power series 
in $t$ with coefficients in $A$.
A {\it \cNcs $($noncommutative symmetric$)$ system} over $A$ 
(see Definition \ref{Main-Def}) 
by definition is a $5$-tuple 
$\Omega\in A[[t]]^{\times 5}$ 
which satisfies the defining equations 
(see Eqs.\,$(\ref{UE-0})$--$(\ref{UE-4})$) 
of the NCSFs (noncommutative symmetric functions) 
first introduced and studied 
in the seminal paper \cite{G-T}. 
When the base algebra 
$K$ is clear in the context,
the ordered pair $(A, \Omega)$ 
is also called a {\it \cNcs system}.
In some sense, a \cNcs system over 
the $K$-algebra $A$ can be viewed 
as a system of analogs in $A$
of the NCSFs defined by 
the same equations.
For more studies on NCSFs, 
see \cite{T}, \cite{NCSF-II}, 
\cite{NCSF-III}, \cite{NCSF-IV}, 
\cite{NCSF-V} and \cite{NCSF-VI}.

One immediate but probably the most 
important example of the \cNcs systems 
is $(\cNsf, \Pi)$ (see Eqs.\,(\ref{lambda(t)})--(\ref{xi(t)})\,) 
formed by the generating functions of 
the NCSFs defined in \cite{G-T}
by Eqs.\,$(\ref{UE-0})$--$(\ref{UE-4})$ 
over the free $K$-algebra $\cNsf$ of NCSFs. 
It serves as the universal \cNcs system 
over all associative $K$-algebras 
(see Theorem \ref{Universal}).  
More precisely, for any \cNcs system $(A, \Omega)$, 
there exists a unique $K$-algebra homomorphism 
$\cS: \cNsf \to A$ such that 
$\cS^{\times 5}(\Pi) = \Omega$ 
(Here we have extended the homomorphism
$\cS$ to $\cS: \cNsf[[t]] \to A[[t]]$ 
by the base extension). 

Note that, it has been shown 
in \cite{G-T}, in the quotient modulo 
the commutator of $\cNsf$, the NCSFs in $\Pi$ 
become the corresponding classical (commutative) 
symmetric functions (\cite{M}). Hence,
the universal \cNcs system for 
commutative $K$-algebras is given by 
the generating functions of 
the corresponding 
classical (commutative) symmetric functions.

One of the main motivations for the introduction 
of \cNcs is as follows (see Subsection \ref{S2.3}
for more discussions).
Note that, as an important topic 
in the theory of symmetric functions, 
the relations or polynomial identities 
among various commutative 
or noncommutative symmetric functions 
have been known explicitly 
(see \cite{M} and \cite{G-T}). 
When a \cNcs system $\Omega$ 
is given over a $K$-algebra $A$, 
by applying the $K$-algebra 
homomorphism $\cS:\cNsf \to A$ 
guaranteed by the universal 
property of the system $(\cNsf, \Pi)$ 
to the identities of the NCSFs 
in $\Pi$, we see the same identities hold
for the elements of $A$ in the \cNcs 
system $\Omega$.
This could be a very effective way 
to obtain identities for certain elements 
of $A$ if we could show they are 
involved in a \cNcs system over $A$.
On the other hand, if the given 
\cNcs system $(A, \Omega)$ has already 
been well-understood, the $K$-algebra 
homomorphism $\cS:\cNsf  \to A$ in turn 
gives a {\it specialization} 
or {\it realization} (\cite{G-T}, \cite{St2})
of NCSFs, which may provide 
some new understanding of NCSFs.
For more studies on the specializations 
of NCSFs, see the references quoted 
above.

This paper is the first of a sequence papers 
on the \cNcs systems over differential operator 
algebras in commutative or noncommutative 
variables (\cite{GTS-II}); the \cNcs systems over
the Grossman-Larson Hopf algebras 
of labeled rooted trees (\cite{GTS-IV});
as well as their connections and applications 
to the inversion problem (\cite{BCW}, \cite{E}) 
and specializations of NCSFs 
(\cite{GTS-III}, \cite{GTS-V}).
In this paper, we first introduce 
the notion {\it \cNcs systems} 
over any associative $K$-algebras.  
We then prove some results on 
the \cNcs systems in general, 
the \cNcs systems over 
bi-algebras or Hopf algebras 
and the universal \cNcs system 
from NCSFs (\cite{G-T}), 
which will be needed 
in the followed papers. 
Finally, we briefly review some of 
the main results that will be proved 
in the followed papers \cite{GTS-II}, 
\cite{GTS-IV} and \cite{GTS-V} 
as some supporting examples 
to the general discussions given 
in the first part of this paper. 

The arrangement of this paper 
is as follows. 
In Subsection \ref{S2.1}, 
we first formulate 
the notion \cNcs systems $(A, \Omega)$ over 
any associative $K$-algebra $A$ 
(see Definition \ref{Main-Def}).
We then show in Lemma \ref{E-Uniq}
the existence and uniqueness of
the solutions in $A[[t]]$ of any one of 
Eqs.\,(\ref{UE-1})--(\ref{UE-4}).
Several direct consequences of Lemma \ref{E-Uniq}
are given in Corollaries 
\ref{Uniqueness-1}--\ref{Uniqueness-4}.
Finally, in
Proposition \ref{bialg-case}, 
we prove a property
of the \cNcs systems over bi-algebras 
or Hopf algebras. 
In Subsection \ref{S2.2}, 
we first recall some NCSFs introduced 
in \cite{G-T} and a graded Hopf algebra 
structure of the space $\cNsf$ of NCSFs.
We then in Theorem \ref{Universal}
show the generating functions
of these NCSFs form the universal
\cNcs system $(\cNsf, \Pi)$ over
all $K$-algebra. Moreover, 
we also give some sufficient conditions 
in Theorem \ref{Universal}
for the algebra homomorphisms guaranteed
by the universal property of 
$(\cNsf, \Pi)$ to be further
homomorphisms of bi-algebras 
and Hopf algebras.
In Subsection \ref{S2.3}, 
we discuss some possible applications of 
the universal properties of  
the \cNcs system $(\cNsf, \Pi)$, 
which are also the main motivations 
for the introduction of 
the \cNcs systems over 
associative $K$-algebras. 

The results above form 
the first part of this paper. 
In the second part, Sections $3$ and $4$, 
we review some of the main results 
that will appear in the sequels 
\cite{GTS-II}, \cite{GTS-IV} and \cite{GTS-V}. 
The main purposes that we include 
Sections $3$ and $4$ in this paper 
are as follows. First, we think 
it is better to provide some concrete 
examples for the general discussions 
of \cNcs systems given in Section \ref{S2}, 
so the paper can be read as 
a more complete introduction 
to the newly defined \cNcs systems. 
Secondly, considering 
length of the whole series of papers,
we hope that discussions in Sections 
\ref{S3} and \ref{S4} 
can also serve as a shorter survey 
or review for some of the main results 
obtained the followed papers 
\cite{GTS-II}, \cite{GTS-IV} and \cite{GTS-V} 
(Precisely speaking, they shold be read 
as an annuouncement since these sequel papers 
for the time being are 
still under submission).


In Section \ref{S3}, 
we discuss the \cNcs systems that will be 
constructed in \cite{GTS-II} over 
differential operator algebras 
in commutative or noncommutative 
free variables. Certain properties 
of the resulting specialization of NCSFs 
by differential operators, which will 
be proved in \cite{GTS-II} 
and \cite{GTS-V}, are also discussed. 
In Section \ref{S4}, for any non-empty $W\subseteq\bN^+$, 
we first recall the Connes-Kreimer Hopf algebra 
$\cH_{CK}^W$ and the Grossman-Larson Hopf algebra 
$\mathcal H^W_{GL}$ of $W$-labeled rooted forests 
and $W$-labeled rooted trees, respectively. 
We then discuss the \cNcs system 
$(\cH_{GL}^W, \Omega_\bT^W)$ 
that will be constructed in \cite{GTS-IV} 
over the Grossman-Larson Hopf algebra 
$\mathcal H^W_{GL}$. Some of 
properties to be given in 
\cite{GTS-IV} and \cite{GTS-V} 
of the resulting specializations 
of NCSFs by $W$-labeled rooted trees
will also be discussed in this section.  
Finally, we briefly explain a connection,
which will be given in \cite{GTS-V}, 
between the \cNcs system 
$(\cH_{GL}^W, \Omega_\bT^W)$ 
with the \cNcs systems discussed in 
Section \ref{S3} over 
differential operator algebras.
Some consequences of this connection 
to the related specializations of NCSFs 
and the inversion problem 
(\cite{BCW}, \cite{E})
will also be discussed.

{\bf Acknowledgment}:
The author would like to thank the referee 
for pointing out some misprints and 
suggesting a connection of \cNcs systems with 
the combinatorial Hopf algebras 
introduced and studied by
Marcelo Aguiar, Nantel Bergeron 
and Frank Sottile in \cite{ABS}.

\renewcommand{\theequation}{\thesection.\arabic{equation}}
\renewcommand{\therema}{\thesection.\arabic{rema}}
\setcounter{equation}{0}
\setcounter{rema}{0}

\section{\bf \cNcs Systems over Associative Algebras} \label{S2}

Let $K$ be any unital commutative $\bQ$-algebra and 
$A$ any unital associative but not necessarily commutative 
$K$-algebra. Let $t$ be a formal central parameter, 
i.e. it commutes with all elements of $A$, and $A[[t]]$ 
the $K$-algebra of formal power series 
in $t$ with coefficients in $A$. 

First let us introduce the following main notion 
of this paper, which is mainly motivated by the 
seminal work \cite{G-T} by I. M. Gelfand; 
D. Krob; A. Lascoux; 
B. Leclerc; V. S. Retakh and
J.-Y. Thibon on NCSFs 
(noncommutative symmetric functions).  

\begin{defi} \label{Main-Def}
For any unital associative $K$-algebra $A$, a $5$-tuple $\Omega=$ 
$( f(t)$, $g(t)$, $d\,(t)$, $h(t)$, $m(t) ) 
\in A[[t]]^{\times 5}$ is said 
to be a {\it \cNcs $($noncommutative symmetric$)$ system} over $A$ 
if the following equations are satisfied.
\begin{align}
&f(0)=1 \label{UE-0}\\
& f(-t)  g(t)=g(t)f (-t)=1, \label{UE-1}   \\
& e^{d\,(t)} = g(t), \label{UE-2} \\
& \frac {d g(t)} {d t}= g(t) h(t), \label{UE-3}\\ 
& \frac {d g(t)}{d t} =  m(t) g(t).\label{UE-4}
\end{align}
\end{defi}

When the base algebra $K$ is clear in the context, we also call 
the ordered pair $(A, \Omega)$ a {\it \cNcs system}. 
Since \cNcs systems often come from generating functions 
of certain elements of $A$ that are under concern, 
the components of $\Omega$ will also be referred as 
the {\it generating functions} of their coefficients.

\subsection{\cNcs Systems in General}\label{S2.1}

In this subsection, we prove some results 
for the \cNcs systems in general and also 
the \cNcs systems over bi-algebras 
and Hopf algebras.

First, let us fix the following 
convention that will be used 
implicitly throughout this paper. 

\vskip3mm

{\bf Convention:}

\vskip2mm
\begin{enumerate}
\item[$(a)$] {\it All $K$-algebras in this paper 
are assumed to be unital associative
$K$-algebras and all $K$-algebra 
homomorphsims are assumed 
to be unit-preserving.}

\item[$(b)$]
{\it For any $K$-algebras $B$ and $A$ and any $K$-linear map 
$\cS: B\to A$, we always extend 
$\mathcal S$ to a linear map, 
which we will still denote by $\mathcal S$, 
from $B[[t]]$ to $A[[t]]$ by the base extension, 
i.e. for any $\sum_{m\geq 0} b_m t^m\in B[[t]]$, 
we set
\begin{align}
\mathcal S(\sum_{m\geq 0} b_m t^m)
=\sum_{m\geq 0} \mathcal S (b_m) t^m.
\end{align}
Furthermore, for any $m \geq 1$, 
we denote by $\cS^{\times m}$ 
the $K$-linear map from 
$B[[t]]^{\times m}$ to 
$A[[t]]^{\times m}$
induced by $\cS: B[[t]]\to A[[t]]$.}
\end{enumerate}
\vskip3mm

Now let $A$ be any unital $K$-algebra and 
$\Omega$ a \cNcs system over $A$ 
as given in Definition \ref{Main-Def}. We define 
five sequences of elements of $A$ by writing
\allowdisplaybreaks{
\begin{align}
f (t)&:=\sum_{m\geq 0} t^m \lambda_m.\label{Def-lambda-m}\\
g(t) &:=\sum_{m\geq 0} t^{m} s_m, \label{Def-s-m} \\
d\,(t)&:=\sum_{m\geq 1} \frac {t^m}m  \phi_m, \label{Def-phi-m}\\
h(t)&:=\sum_{m\geq 1} t^{m-1} \psi_m, \label{Def-psi-m}\\
m(t)&:=\sum_{m\geq 1} t^{m-1} \xi_m. \label{Def-xi-m}
\end{align}

We will also denote each sequence of 
the elements of $A$ defined above 
by the corresponding letter without sub-index. 
For example, $\lambda$ denotes the sequence 
$\{\lambda_m\,|\, m\geq 0\}$ defined in 
Eq.\,(\ref{Def-lambda-m}) and $\xi$ denotes 
the sequence $\{\xi_m\,|\, m\geq 1 \}$ 
defined in Eq.\,(\ref{Def-xi-m}), etc.

Next, let us start with the following  
simple but useful properties of \cNcs systems.

\begin{lemma}\label{Hom-NS}
Let $A$ and $B$ be any $K$-algebras and $\cS \! :\! B\to A$ 
an algebra homomorphism. Let $\tilde \Omega$ be a \cNcs system over $B$.
Then $\Omega:=\cS^{\times 5} (\tilde \Omega)$ is a \cNcs system over $A$.
\end{lemma}

\pf Note that, by our convention above,  
$\cS: B[[t]] \to A[[t]]$ is also a unital $K$-algebra homomorphism 
and commutes with the linear operator $\frac{d}{d\, t}$. 
With these observations, it is easy to see that, 
the components of the $5$-tuple 
$\Omega:= \cS^{\times 5} (\tilde \Omega)$ also satisfy 
Eqs.\,(\ref{UE-0})--(\ref{UE-4}). 
Hence the lemma follows.
\epfv

\begin{lemma}\label{flip}
Let $(A, \Omega)$ be a \cNcs system as fixed above. Define
the $5$-tuple $\Omega^\tau$ to be 
\begin{align}\label{flip-eq}
\Omega^\tau := (\, g(-t), \, f(-t), \, -d\,(t),\,  -m(t), \, -h(t) )\,. 
\end{align}
Then $\Omega^\tau$ is also a \cNcs system over $A$.
\end{lemma}

We call the \cNcs system $\Omega^\tau$ the {\it flip} of $\Omega$.

\pf First, 
for convenience, we also write $\Omega^\tau$ as
\begin{align*}
\Omega^\tau=:( \, \tilde f(t),\, \tilde g(t),\, 
\tilde d\,(t), \, \tilde h(t), \, \wtilde m(t)\,)\,, 
\end{align*}
i.e. we set $\tilde f(t):=g(-t)$; $\tilde g(t):=f(-t)$; etc.

By setting $t=0$ in  
Eq.\,(\ref{UE-1}) for $\Omega$, 
we see that $\tilde f(0)=g(0)=1$. 
So we get Eq.\,(\ref{UE-0}) for $\Omega^\tau$. 
By rewriting Eq.\,(\ref{UE-1}) for $\Omega$ in terms  
of $\tilde f(t)$ and $\tilde g(t)$, we see 
Eq.\,(\ref{UE-1}) also holds
for $\wtilde \Omega$. 
To show Eq.\,(\ref{UE-2}) for $\Omega^\tau$, 
we consider 
\begin{align*}
e^{\tilde d(t)}&=e^{-d(t)}\\
&=g(t)^{-1}\\
&=f(-t)\\
&=\tilde g(t).
\end{align*}
Hence we get Eq.\,(\ref{UE-2}) 
for $\Omega^\tau$. 

Now consider Eqs.\,(\ref{UE-3}) and (\ref{UE-4})
for $\Omega^\tau$. 
Note first that, 
by applying 
$\frac{d}{d\, t}$ to Eq.\,(\ref{UE-1}) 
for $\Omega$, we have
\begin{align*}
0&=\frac{df(-t)}{d\, t}g(t)+ f(-t)\frac{d g(t)}{d\, t}\, ,\\
0&=\frac{d g(t)}{d\, t}f(-t)+
g(t)\frac{df(-t)}{d\, t}.
\end{align*}
Replacing $\frac{d g(t)}{d\, t}$  
by $m(t)g(t)$ and $g(t)h(t)$ respectively
in the last two equations above
and then solving $\frac{df(-t)}{d\, t}$, 
we get
\begin{align*}
\frac{d  f(-t)}{d\, t} & =-f(-t)m(t)\, ,\\ 
\frac{d f(-t)}{d\, t} & =-h(t)f(-t)\, .
\end{align*}

Rewriting the last two equations above in terms of 
$\tilde g(t)$, $\tilde h(t)$ and $\wtilde m(t)$, 
we get Eqs.\,(\ref{UE-3}) 
and $(\ref{UE-4})$ for 
$\Omega^\tau$.
\epfv

Let $A_i$ $(i=1, 2)$ be any $K$-algebras and $\Omega_i$ 
a \cNcs system over $A_i$. Write $\Omega_i$ as
\begin{align*}
\Omega_i=(f_i(t), \, g_i(t),\, d_i(t),\, 
h_i(t),\, m_i(t))
\end{align*}
and set 
\allowdisplaybreaks{
\begin{align}
f_3(t)   &:=f_1(t) \otimes f_2(t) \, ,  \nno \\
g_3(t)   &:=g_1(t) \otimes g_2(t)\, ,\nno \\
d_3(t)   &:= d_1(t) \otimes 1 +  1 \otimes d_2(t)\, ,\nno \\
h_3(t)   &:= h_1(t) \otimes 1 +  1 \otimes h_2(t)\, ,\nno \\
m_3(t)   &:= m_1(t) \otimes 1 +  1\otimes m_2(t)\, \nno \\
\Omega_1 \otimes_K \Omega_2 &:=(f_3(t),  \, g_3(t),\, d_3(t)\, , \label{tensor-O} 
h_3(t),\, m_3(t))\, .
\end{align} }

We call $\Omega_1 \otimes_K \Omega_2 $ the 
{\it tensor product} (over $K$) 
of the \cNcs systems 
of $\Omega_1$ and $\Omega_2$.
Then we have the following proposition.

\begin{propo}\label{Tensor-NS}
Let $(A_i, \Omega_i)$ $(i=1, 2)$ be
\cNcs systems.  Then, the tensor product 
$\Omega_1 \otimes_K \Omega_2$  
 of the \cNcs systems 
of $\Omega_1$ and $\Omega_2$
forms a \cNcs system over the $K$-algebra 
$A_1\otimes_K A_2$.
\end{propo}

\pf The proof is straightforward, and is just 
to check the components 
of $\Omega_1 \otimes_K \Omega_2$ satisfy 
Eqs.\,(\ref{UE-0})--(\ref{UE-4}).
First, it is easy to see that Eqs.\,(\ref{UE-0}) and 
(\ref{UE-1}) are satisfied. To show Eq.\,(\ref{UE-2}), note that, 
$d_1(t)\otimes 1$ and $1\otimes d_2(t)$ as elements 
of $(A_1\otimes_K A_2) [[t]]$ commute with each other. 
By this fact, we have
\begin{align*}
e^{d_3(t)}&= e^{d_1(t) \otimes 1 +  1 \otimes d_2(t)} \\
&=e^{d_1(t) \otimes 1} e^{1 \otimes d_2(t)} \\
\intertext{Applying Eq.\,(\ref{UE-2}) for $g_1(t)$ and $g_2(t)$:}
&=(g_1(t)\otimes 1)(1 \otimes g_2(t))\\
&=g_1(t) \otimes g_2(t)\\
&=g_3(t).
\end{align*}

Next, let us show Eq.\,(\ref{UE-3}) as follows.  
\begin{align*}
\frac {d\, g_3(t) }{d\, t} &= \frac {d}{d \, t} \, ( g_1(t) \otimes g_2(t) )\\
&=  \frac {d\, g_1 (t) }{d\, t}  \otimes g_2(t) 
+ g_1 (t)\otimes  \frac {d\, g_2(t) }{d\, t} \\
\intertext{Applying Eq.\,(\ref{UE-3}) for $g_1(t)$ and $g_2(t)$:}
&=(g_1(t) h_1(t))\otimes g_2(t) + g_1(t) \otimes (g_2(t) h_2(t)) \\
&=(g_1(t) \otimes g_2(t))( h_1(t) \otimes 1 +  1 \otimes h_2(t) )\\
&=g_3(t) h_3(t).
\end{align*}

Hence, we get Eq.\,(\ref{UE-3}).
Eq.\,(\ref{UE-4}) can be proved similarly.
\epfv

Next, let us prove the existence and uniqueness 
of solutions for any one of 
Eqs.\,$(\ref{UE-1})$--$(\ref{UE-4})$.

\begin{lemma}\label{E-Uniq}
Let $A$ be any associative $K$-algebra. Then, 
for any one of Eqs.\,$(\ref{UE-1})$--$(\ref{UE-4})$,
if one generating function in the equation 
fixed to be an element of $A[[t]]$ 
with $f(0)=1$, $g(0)=1$ or $d(0)=0$ 
if $f(t)$, $g(t)$ or $d(t)$ is the one fixed,
the equation has one and only one solution in $A[[t]]$
for the other generating function.
\end{lemma}

\pf First, the lemma is obvious for 
Eqs.\,$(\ref{UE-1})$ and $(\ref{UE-2})$.
To see it is also true for Eq.\,$(\ref{UE-3})$, 
we write $g(t)$ and $h(t)$ as in 
Eqs.\,(\ref{Def-s-m}) and (\ref{Def-psi-m}), respectively.
Using the fact $s_0=1$ and Eq.\,(\ref{UE-3}), we have 
\begin{align*} 
\sum_{m\geq 1} m s_m t^{m-1} = 
\lp 1 + \sum_{m\geq 1} s_m t^{m}\rp
\lp \sum_{m\geq 1} \psi_m t^{m-1} \rp.
\end{align*}

Comparing the coefficients of $t^{m-1}$ $(m\geq 2)$ 
in the equation above, we get
\begin{align} 
\psi_1&=s_1,\label{Uniqueness-pe1}\\ 
m s_m &=\psi_m+ \sum_{\substack{k+l=m, \\
k, l\geq 1 }}  s_k \psi_l, \label{Uniqueness-pe2}
\end{align}
for any $m\geq 2$. 

Then, by Gauss' elimination method, 
it is easy to see 
that, if one of the sequences $\{s_m\,|\, m\geq 1\}$ and 
$\{\psi_m\,|\, m\geq 1\}$ is given, 
the other can always be obtained in a unique way. 
Hence the lemma is true for Eq.\,$(\ref{UE-3})$.
Similarly, we can prove the lemma 
for Eq.\,$(\ref{UE-4})$.
\epfv

From Lemma \ref{E-Uniq} and its proof, 
it is easy to see that 
we have the following three corollaries.

\begin{corol}\label{Uniqueness-1}
For any $K$-algebra $A$ and $c(t)\in A[[t]]$, we have

$(a)$ If $c(0)=1$, then, for any $i=1$ or $2$, 
there exists a unique \cNcs system $\Omega$ 
over $A$ with $c(t)$ as the $i^{th}$ component.

$(b)$ If $c(0)=0$, then, there exists 
a unique \cNcs system $\Omega$ 
over $A$ with $c(t)$ 
as the $3rd$ component.

$(c)$ For any $i=4$ or  $5$, there exists a 
unique \cNcs system $\Omega$ 
over $A$ with $c(t)$ as the $i^{th}$ 
component.
\end{corol}

\begin{corol}\label{Uniqueness-2}
Let $(A, \Omega)$ be a \cNcs system. 
Then, any component of $\Omega$ 
completely determines the others. 
In other words, 
if two \cNcs systems over $A$
have a same component at a same location,
then these two systems are completely same.
\end{corol}

\begin{corol}\label{Uniqueness-3}
Let $(A, \Omega)$ be a \cNcs system. 
For any sequence $w=\{w_m\,|\, m\geq 1 \}$ 
of elements of $A$,
we denote by $K \langle w \rangle$ 
the unital subalgebra of $A$ generated by $w_m$'s. 
Then, for any sequence $w=s$, $\phi$, $\psi$ or $\xi$, 
we have
\begin{align}
K\langle w \rangle=K \langle \lambda \rangle \, .
\end{align}
\end{corol}

\begin{corol}\label{Uniqueness-4}
Let $(A, \Omega)$ and $(B, \tilde \Omega)$ be \cNcs systems 
and $\mathcal S: B \to A$ a $K$-algebra homomorphism.
Suppose that, for some  $1\leq j\leq 5$, $\cS$ maps 
the $j^{th}$ component of $\tilde \Omega$ 
to the $j^{th}$ component of $\Omega$. 
Then,  $\mathcal S^{\times 5} (\tilde \Omega)=\Omega$.
\end{corol}

\pf First, by Lemma \ref{Hom-NS}, 
we see that $\cS^{\times 5}(\tilde \Omega)$ 
also is a \cNcs system over $A$.
Since the \cNcs systems $\cS^{\times 5}(\tilde\Omega) $ 
and $\Omega$ over $A$ have same $j^{th}$ 
component, by Corollary \ref{Uniqueness-2}, 
we have $\mathcal S^{\times 5} (\tilde\Omega)=\Omega$.
\epfv

\begin{propo}\label{tau-flip}
Let $(A, \Omega)$ a \cNcs system as fixed 
in Definition \ref{Main-Def} and 
$\tau: A \to A$ a $K$-algebra homomorphism
such that $\tau (\phi_m)=-\phi_m$ for any $m\geq 1$.
Then, we have  
\begin{align}
\tau (\lambda_m)&=(-1)^m s_m   \, ,\label{tau-flip-e1} \\
\tau (s_m)&=(-1)^m \lambda_m \, ,\label{tau-flip-e2} \\
\tau (\psi_m)&=- \xi_m \, ,\label{tau-flip-e3}  \\
\tau (\xi_m)&=-\psi_m\, .\label{tau-flip-e4} 
\end{align}
\end{propo}

\pf Let  $\tau^{\times 5}(\Omega)$ and $\Omega^\tau$ be 
respectively the image of $\Omega$ 
under the homomorphism $\tau$ 
and the flip of $\Omega$ defined in 
Eq.\,(\ref{flip-eq}). 
By lemma \ref{Hom-NS} and \ref{flip}, we know both 
$\tau^{\times 5}(\Omega)$ and $\Omega^\tau$ 
are \cNcs systems over $A$. While on the other hand, 
by the conditions in the proposition, we know 
$\tau^{\times 5}(\Omega)$ and $\Omega^\tau$ have the 
same third component $-d(t)$. 
So, by Corollary \ref{Uniqueness-2}, 
we have $\tau^{\times 5} (\Omega)=\Omega^\tau$ 
from which Eqs.\,(\ref{tau-flip-e1})--(\ref{tau-flip-e4}) 
follow directly. 
\epfv

Next we consider \cNcs systems over some special algebras $A$.

\begin{propo}\label{Comm-Case}
Let $(A, \Omega)$ be any \cNcs system 
over a commutative $K$-algebra $A$. 
Then we have 
\begin{align}
m(t)=h(t)=d'(t), \label{Comm-Case-e1}
\end{align}
where $d'(t)$ denotes the first 
derivative of $d(t)$ over $t$.

In particular, the system of Eqs.\,$(\ref{UE-0})$--$(\ref{UE-4})$
is reduced to:
\begin{align*}
\begin{cases}
 &f(0)=1 \, , \\
& f(-t)g(t)=1\, ,   \\
&  \frac {d g(t)} {d t}= g(t) h(t)\, , \quad \text{or equivalently,} \quad 
e^{d\,(t)} = g(t)\, .
\end{cases}
\end{align*}
\end{propo}

\pf Since $A$ is commutative, so is $A[[t]]$. In this case, 
Eqs.\,(\ref{UE-3}) and (\ref{UE-4}) become same.
Applying $\fr{d}{d t}$ 
to Eq.\,(\ref{UE-2}) and, by the chain rule, we have
\begin{align}
\fr{d\, g(t)}{d\, t}=\fr{d}{d t} e^{d(t)}=e^{d(t)}d'(t)=g(t)d'(t). \nno
\end{align}

From the observations above, we see that $h(t)$, $m(t)$ and
$d'(t)$ are all solutions of Eq.\,(\ref{UE-3}) with (same) 
$g(t)$. Therefore, by 
Lemma \ref{E-Uniq}, we have Eq.\,(\ref{Comm-Case-e1}) 
and the proposition follows.
\epfv

Next, we consider \cNcs systems over $K$-bialgebras. 
First we need recall the following notions 
(see \cite{Abe}, \cite{Knu} 
and \cite{Mon} for more details). 
Let $A$ be a $K$-bialgebra with 
the co-product denoted by $\Delta: A\to A\otimes A$.
An element $x\in A$ is  
{\it primitive} if $\Delta (x)=1\otimes x+x\otimes 1$.
$x\in A$ is a {\it group-like} element if $\Delta (x)=x\otimes x$.
A sequence $\{a_m\,|\, m\geq 0\}$ of elements of A is said 
to be {\it a sequence of divided powers} if, 
for any $m\geq 0$, we have
\begin{align}\label{Divided}
\Delta a_m=\sum_{\substack{k+l=m \\k, l\geq 0}} a_k\otimes a_l.
\end{align}

Now let $t$ be a central parameter as before, 
we extend the counit $\epsilon$ of $A$ 
to $\epsilon:A[[t]] \to k$ by setting $\epsilon(t)=0$
and the co-product $\Delta$ of $A$
to $\Delta:A[[t]]\to A[[t]]\otimes_K A[[t]]$ 
by the base extension. Then, 
with the extended counit and co-product, 
$A[[t]]$ is also a $K$-bialgebra. 
With this $K$-bialgebra structure fixed on $A[[t]]$, 
for any sequence $\{a_m \,|\, m\geq 0\}$ 
of elements of $A$,
it is easy to check that the following facts: 
\begin{enumerate}
\item[$\bullet$]
the sequence $\{a_m  \,|\, m\geq 0\}$ 
is a sequence of divided powers of $A$ 
iff its generating function 
$a(t):=\sum_{m\geq 0} a_m t^m$ 
is a group-like element of $A[[t]]$;
\item[$\bullet$]
all elements $a_m$ $(m \geq 0)$ 
are primitive in $A$ 
iff the generating function
$a(t)$ is a primitive element 
of $A[[t]]$.
\end{enumerate}

\begin{propo}\label{bialg-case}
Let $(A, \Omega)$ be a \cNcs system. 
Suppose $A$ is further a $K$-bialgebra. 
Then the following statements are equivalent. 

$(1)$ $\{\lambda_m \,|\, m\geq 0\}$ is  a sequence of divided powers of $A$.

$(2)$ $\{s_m \,|\, m\geq 0\}$ is  a sequence of divided powers of $A$.

$(3)$ One $(\text{hence also all})$ of \,$d(t)$, $h(t)$ and $m(t)$ 
is primitive in $A[[t]]$.
\end{propo}

Note that, the statement $(3)$ is same as saying that,
the sequence $\{\phi_m\,|\, m\geq 1\}$,  $\{\psi_m\,|\, m\geq 1\}$
or $\{\xi_m\,|\, m\geq 1\}$ is a sequence of 
the primitive elements of $A$.

\pf By the discussion before the proposition, it will be enough 
to show that the following equations are equivalent to each other.
\allowdisplaybreaks{
\begin{align}
\Delta f(t)&=f(t)\otimes f(t), \label{bialg-case-pe1}\\
\Delta g(t)&=g(t)\otimes g(t),\label{bialg-case-pe2}\\
\Delta d(t)&=d(t)\otimes 1+1\otimes d(t),\label{bialg-case-pe3}\\
\Delta h(t)&=h(t)\otimes 1+1\otimes h(t),\label{bialg-case-pe4}\\
\Delta m(t)&=m(t)\otimes 1+1\otimes m(t).\label{bialg-case-pe5}
\end{align} }

First, we identify the $K$-algebras
$A[[t]]\otimes_K A[[t]]$ with $(A \otimes_K A)[[t]]$ 
in the standard way. Then, both sides of 
Eqs.\,(\ref{bialg-case-pe1})--(\ref{bialg-case-pe5})
can be viewed as elements of the $K$-algebra 
$(A \otimes_K A)[[t]]$. 

Secondly, note that, 
the $5$-tuple of $(A \otimes_K A)[[t]]$ formed by the LHS's 
of Eqs.\,(\ref{bialg-case-pe1})--(\ref{bialg-case-pe5})
in the same order as the equations displayed above
is the image $\Delta^{\times 5} (\Omega)$
in $\lp (A\otimes A)[[t]]\rp^{\times 5}$
of the \cNcs system $\Omega$ over $A$
under the $K$-algebra homomorphism 
$\Delta^{\times 5} : A^{\times 5} \to (A\otimes_K A)^{\times 5}$; 
while the $5$-tuple of $(A \otimes_K A)[[t]]$ on the RHS's is 
the tensor product $\Omega\otimes_K \Omega$ of 
the \cNcs system $\Omega$ with itself.
Then, by Lemma \ref{Hom-NS}, we know $\Delta \Omega$
is a \cNcs system over the $K$-algebra $(A \otimes_K A)[[t]]$, 
and,  by Proposition \ref{Tensor-NS},
$\Omega\otimes_K \Omega$ is also a \cNcs system over 
$(A \otimes_K A)[[t]]$. Also note that, one of 
Eqs.\,(\ref{bialg-case-pe1})--(\ref{bialg-case-pe5})
holds iff the \cNcs systems $\Delta \Omega$ and 
$\Omega\otimes \Omega$ have a same component 
at a same location. Hence, by Corollary \ref{Uniqueness-2}, 
Eqs.\,(\ref{bialg-case-pe1})--(\ref{bialg-case-pe5})
are equivalent to each other 
and the proposition follows.
\epfv

\subsection{The Universal \cNcs System from
Noncommutative Symmetric Functions}\label{S2.2}

In this subsection, we first recall the definitions 
of some NCSFs (noncommutative symmetric functions) 
first introduced and studied in \cite{G-T}, 
whose generating functions form  
a \cNcs system $\Pi$ over the free associative 
algebra $\cNsf$ generated by 
an alphabet $\{\Lambda_m\,|\, m\geq 1\}$.
We then show in Theorem \ref{Universal} that 
the \cNcs system $(\cNsf, \Pi)$ from NCSFs 
is actually the universal \cNcs system over 
all associative $K$-algebras. When $A$ is further 
a $K$-bialgebra (resp.\,\,Hopf algebra), 
some sufficient conditions for 
the algebra homomorphism 
$\cS: \cNsf\to A$ to be a bialgebra 
(resp.\,\,Hopf algebra) homomorphism 
are also given in 
Theorem \ref{Universal}.

First, let $K$ be a unital commutative $\bQ$-algebra as before and 
$\Lambda=\{ \Lambda_m\,|\, m\geq 1\}$ 
be an alphabet, i.e. a sequence of noncommutative 
free variables. 
For convenience, we also set $\Lambda_0=1$.
Let $\cNsf$ or $K\langle \Lambda \rangle$ 
be the free associative algebra 
generated by $\Lambda$ over $K$. 
We denote by
$\lambda (t)$ the generating function of 
$\Lambda_m$ $(m\geq 0)$, i.e. we set

\begin{align}\label{lambda(t)}
\lambda (t):= \sum_{m\geq 0} t^m 
\Lambda_m =1+\sum_{k\geq 1} t^m \Lambda_m.
\end{align}

In the theory of NCSFs, 
$\Lambda_m$ $(m\geq 0)$ is 
the noncommutative analog 
of the $m^{th}$ classical (commutative) 
elementary symmetric function 
and is called the {\it $m^{th}$ 
$(\text{noncommutative})$ 
elementary symmetric function.}

To define some other NCSFs, we consider 
Eqs.\,$(\ref{UE-1})$--$(\ref{UE-4})$ 
over the free $K$-algebra $\cNsf$
with $f(t)=\lambda(t)$. The 
solutions for $g(t)$, $d\,(t)$, 
$h(t)$, $m(t)$ exist and are unique 
(see Corollary \ref{Uniqueness-1}, for example), 
whose coefficients will be the NCSFs 
that we are going to define.
Following the notation in \cite{G-T}, 
we denote the resulting 
$5$-tuple by 
\begin{align}\label{Def-Pi}
\Pi=(\lambda(t),\, \sigma(t),\, \Phi(t),\, \psi(t),\, \xi(t))
\end{align}
and write the last 
four generating functions of 
$\Pi$ explicitly as follows.

\allowdisplaybreaks{
\begin{align}
\sigma (t)&=\sum_{m\geq 0} t^m S_m,  \label{sigma(t)} \\
\Phi (t)&=\sum_{m\geq 1} t^m \frac{\Phi_m}m  \label{Phi(t)}\\
\psi (t)&=\sum_{m\geq 1} t^{m-1} \Psi_m, \label{psi(t)}\\
\xi (t)&=\sum_{m\geq 1} t^{m-1} \Xi_m.\label{xi(t)}
\end{align}}

Note that, in terms of the terminology 
in the previous subsection, the $5$-tuple $\Pi$ 
defined above is the unique \cNcs system 
with $f(t)=\lambda(t)$ in 
Eq.\,(\ref{lambda(t)}) 
over the free $K$-algebra 
$\cNsf$.

Following \cite{G-T}, we call $\S_m$ $(m\geq 1)$ the {\it $m^{th}$ 
complete homogeneous symmetric function}, and 
$\Psi_m$ and $\Xi_m$ $(m\geq 1)$ respectively the
{\it $m^{th}$ power sum symmetric 
function of the first and second kind}. 
Note that, $\Xi_m$ $(m\geq 1)$ were denoted by $\Psi_m^*$ 
in \cite{G-T}. Due to Proposition \ref{omega-Lambda} 
below,  the NCSFs $\Xi_m$ $(m\geq 1)$
do not play an important role in the NCSF theory 
(see the comments in page $234$ in \cite{G-T}). 
But, in the context of some other problems, 
relations of $\Xi_m$'s with other NCSFs, 
especially, with $\Psi_m$'s, 
are also important. For example, 
this is indeed the case in \cite{GTS-III} 
where connections of NCSFs with the inversion problem 
are concerned. So we here refer $\Xi_m \in \cNsf$ 
$(m\geq 1)$ as the {\it $m^{th}$ $(\text{noncommutative})$ 
power sum symmetric function of the third kind}.

The following two propositions proved in \cite{G-T} 
and \cite{NCSF-II} will be very useful 
for our later arguments.

\begin{propo}\label{bases}
For any unital commutative $\bQ$-algebra $K$,
the free algebra $\cNsf$ is freely generated
by any one of the families of the NCSFs 
defined above.
\end{propo}


\begin{propo}\label{omega-Lambda}
Let $\omega_\Lambda$ be the anti-involution of 
$\cNsf$ 
which fixes $\Lambda_m$ $(m\geq 1)$.
Then, for any $m\geq 1$, we have
\begin{align}
\omega_\Lambda (S_m)&=S_m,  \label{omega-Lambda-e1}\\
\omega_\Lambda (\Phi_m)&=\Phi_m,\label{omega-Lambda-e2} \\
\omega_\Lambda (\Psi_m)&=\Xi_m. \label{omega-Lambda-e3}
\end{align}
\end{propo}

As shown in \cite{G-T}, 
the connections between the NCSFs 
and the classical (commutative) 
symmetric functions (\cite{M}), 
are as follows. 
Let $X=\{ X_m\,|\, m\geq 1 \}$ 
be another alphabet and $K \langle X \rangle$ 
the free associative algebra 
generated by $X$ over $K$. We can view 
$\cNsf$ as a subalgebra of $K \langle X \rangle$ 
by setting, for any $m\geq 1$, 
\begin{align}
\Lambda_m=\sum_{i_1<i_2<\cdots<i_m}X_{i_1}X_{i_2}\cdots X_{i_m}.
\end{align}
Let $\mathcal I$ be the two-sided ideal generated 
by the commutators of $X_m$'s and $x$ 
the image of $X$ in the quotient algebra 
modulo $\mathcal I$. Then, 
in the quotient algebra $K[x]$,
$\Lambda_m$ and $S_m$ $(m\geq 1)$ become 
the {\it $m^{th}$ elementary symmetric 
function} and the {\it $m^{th}$ 
complete elementary symmetric 
function}, respectively; 
while $\Phi_m$, $\Psi_m$ and $\Xi_m$ $(m\geq 1)$ 
all become the {\it $m^{th}$ 
power sum symmetric 
function}. See \cite{M} for more studies on the 
classical symmetric functions above.

Next, let us recall the following graded
$K$-Hopf algebra structure 
of $\cNsf$. It has been shown in 
\cite{G-T} that $\cNsf$ is the universal enveloping algebra 
of the free Lie algebra generated 
by $\Psi_m$ $(m\geq 1)$. Hence, it has a Hopf  
$K$-algebra structure as all other 
universal enveloping algebras 
of Lie algebras do. Its co-unit $\epsilon:\cNsf \to K$,
 co-product $\Delta$ and 
 antipode $S$ are uniquely determined by 
\begin{align}
\epsilon (\Psi_m)&=0, \label{counit} \\
\Delta (\Psi_m) &=1\otimes \Psi_m 
+\Psi_m\otimes 1, \label{coprod}\\
S(\Psi_m) & =-\Psi_m,\label{antipode}
\end{align}
for any $m\geq 1$. 

Furthermore, we define the {\it weight} for NCSFs 
by setting the weight of 
any monomial $\Lambda_{m_1}^{i_1} 
\Lambda_{m_2}^{i_2} \cdots \Lambda_{m_k}^{i_k}$
to be $\sum_{j=1}^k i_j m_j$. 
For any $m\geq 0$, we denote by $\cNsf_{[m]}$ 
the vector subspace of $\cNsf$ spanned 
by the monomials of $\Lambda$ 
of weight $m$. Then it is easy to see that 
\begin{align}\label{Grading-cNsf}
\cNsf=\bigoplus_{m\geq 0} \cNsf_{[m]}, 
\end{align}
which provides a grading for $\cNsf$. 

Note that, it has been shown in \cite{G-T}, 
for any $m\geq 1$, the NCSFs 
$S_m, \Phi_m, \Psi_m \in  \cNsf_{[m]}$. 
By Proposition \ref{omega-Lambda}, 
this is also true for the NCSFs $\Xi_m$'s.
By the facts above and 
Eqs.\,(\ref{counit})--(\ref{antipode}), 
it is also easy to check that, 
with the grading given in Eq.\,(\ref{Grading-cNsf}), 
$\cNsf$ forms a graded $K$-Hopf algebra. 
Its graded dual is given 
by the space $\cQf$ of quasi-symmetric functions, 
which were first introduced by I. Gessel \cite{Ge} 
(see also \cite{MR} and \cite{St2} 
for more discussions).

Now we come back to our discussions on the \cNcs systems. 
Note that,  we have seen that $(\cNsf, \Pi)$ by definition 
forms a \cNcs system. More importantly, 
we have the following theorem on the \cNcs system 
$(\cNsf, \Pi)$. 

\begin{theo}\label{Universal}
Let $A$ be a $K$-algebra and $\Omega$ 
a \cNcs system over $A$. Then, 

$(a)$ There exists a unique $K$-algebra homomorphism 
$\cS: \cNsf\to A$ such that 
$\cS^{\times 5} (\Pi)=\Omega$.

$(b)$ If $A$ is further a $K$-bialgebra $($resp.\,\,$K$-Hopf algebra$)$ 
and one of the equivalent statements in Proposition \ref{bialg-case} 
holds for the \cNcs system $\Omega$, then $\cS: \cNsf\to A$ is also 
a homomorphism of $K$-bialgebras $($resp.\,\,$K$-Hopf algebras$)$.
\end{theo}

\pf  $(a)$ Let $\Omega$ be given as in Definition \ref{Main-Def} 
and its components given as in 
Eqs.\,(\ref{Def-lambda-m})--(\ref{Def-xi-m}).
Let $\cS: \cNsf\to A$
to be the unique $K$-algebra 
homomorphism such that, for any $m\geq 1$,
\begin{align}\label{Universal-pe1}
\mathcal S (\Lambda_m)=\lambda_m.
\end{align}

Since $\cNsf$ is freely generated by
$\Lambda_m$ $(m\geq 1)$ 
as a $K$-algebra, 
so $\cS$ is well defined.
Then, by Corollary \ref{Uniqueness-4}, 
we have $\cS^{\times 5} (\Pi)=\Omega$.
The uniqueness of $\cS$ follows 
from the requirement $\cS(\lambda(t))=f(t)$, 
which is same as Eq.\,$(\ref{Universal-pe1})$, 
and again the fact that $\cNsf$ is freely 
generated by $\Lambda_m$ $(m\geq 1)$.

$(b)$ Let $\epsilon$ and $\epsilon_A$ 
denote the co-units  
of $\cNsf$ and $A$, respectively, and 
$\Delta$ denote the co-products 
of both $\cNsf$ and $A$. Then  
we need show the following two equations.
\begin{align}
(\cS \otimes \cS )\circ \Delta 
&= \Delta \circ  \cS 
\label{Universal-pe2} \\
\epsilon_A \circ \cS&=\epsilon. \label{Universal-pe3}
\end{align}

First, by the second condition in $(b)$ and 
Proposition \ref{bialg-case}, we may assume 
$\psi_m$  $(m\geq 1)$ are all 
primitive elements of $A$.
By Eq.(\ref{coprod}), we know that 
$\Psi_m$  $(m\geq 1)$ are all 
primitive elements of $\cNsf$.  
Secondly, note that, 
all the maps involved in 
Eqs.\,$(\ref{Universal-pe2})$ and $(\ref{Universal-pe3})$ 
are $K$-algebra homomorphisms, and,  
by Proposition \ref{bases},
$\cNsf$ is freely generated by 
$\Psi_m$ $(m\geq 1)$ as a $K$-algebra. Therefore, 
to show Eqs.\,$(\ref{Universal-pe2})$ and $(\ref{Universal-pe3})$, 
it will be enough to show that both sides of the equations 
have same values at $\Phi_m$ $(m\geq 1)$. 

With the observations above, for any $m\geq 1$, we consider
\allowdisplaybreaks{
\begin{align*}
(\cS \otimes \cS) (\Delta \Psi_m)
& =(\cS \otimes \cS) (\Psi_m\otimes 1+1\otimes \Psi_m)\\
& =\cS (\Psi_m) \otimes 1 + 1\otimes \cS (\Psi_m)\\
& =\psi_m\otimes 1+1\otimes \psi_m \\
&=\Delta \psi_m.
\end{align*}}
Hence, we have Eqs.\,$(\ref{Universal-pe2})$.

To show Eq.\,$(\ref{Universal-pe3})$, first, 
by Theorem $2.1.3$ in \cite{Abe}, we know 
that the counit of any $K$-bialgebra $B$ 
maps any primitive element $y\in B$ to zero.
Therefore, for any $m\geq 1$, we have 
$\epsilon(\psi_m)=0$ and 
$\epsilon_A (\cS(\Psi_m))
=\epsilon_A (\psi_m)=0$. 
Hence we have Eq.\,$(\ref{Universal-pe3})$.

Finally, let us consider the case 
that $A$ is further a $K$-Hopf algebra.
We need show that $\cS$ in this case 
also commutes with the antipodes 
of $\klam$ and $A$, 
i.e.
\begin{align}
S \circ \cS =\cS \circ S, \label{Universal-pe4}
\end{align}
where both the antipodes of $\klam$ and $A$ 
are denoted by $S$.

First, since both antipodes $S$  
are anti-homomorphisms of $K$-algebras and 
$\cS:\cNsf \to A$ is a $K$-algebra 
homomorphism, it will be enough 
to show that both sides of 
Eq.\,(\ref{Universal-pe4})
have same values 
at $\Psi_m$ $(m\geq 1)$.
Secondly,  it is well-known 
(see \cite{Abe}, for example)
and also easy to check 
that the antipode $S$ of any Hopf algebra $A$ maps
any primitive element $x\in A$ to $-x$.
By the observations above, we have
\begin{align*} 
S \circ \cS \,(\Psi_m)&= S(\psi_m)=-\psi_m,\\
\cS \circ S \, (\Psi_m)&=\cS (-\Psi_m)=-\psi_m.
\end{align*}
Hence, Eq.\,(\ref{Universal-pe4}) holds.
\epfv

\begin{rmk}\label{Universal-Comm}
By taking the quotient over 
the two-sided ideal generated 
by the commutators 
of $\Lambda_m$'s, or 
applying the similar arguments as in the 
proof of Theorem \ref{Universal},
it is easy to see that, 
over the category of commutative $K$-algebras, 
the universal \cNcs system 
is given by the generating functions of
the corresponding classical $($commutative$)$
symmetric functions $($\cite{M}$)$.
\end{rmk}

\begin{rmk}\label{CHA}
Following the referee's suggestion, 
we would like to point out 
a connection of the universal homomorphism 
$\cS: \cNsf \to A$ with 
the universal homomorphism from 
the Hopf algebra $\cQf$ of 
quasi-symmetric functions to 
{\it combinatorial Hopf algebras} 
 introduced and studied in \cite{ABS}.

Suppose that $A$ is a graded 
and connected, and 
one of the statements of 
Proposition \ref{bialg-case} holds, 
say statement $(b)$. 
Furthermore assume in this case 
that the elements $s_m$ $(m\geq 1)$ 
are homogeneous and with grading $m$.
Denote by $A^*$ the graded dual Hopf algebra 
of $A$.  Then, by taking duals, we get 
a homomorphism 
$\cS^*: A^* \to \cQf$
of K-Hopf algebras.
Let $\zeta$ be the linear functional 
of $A^*$ induced by
the sequence $\{s_m\,|\, m\geq 0\}$, i.e.
$\zeta |_{A_m^*}$ $(m\geq 0)$ 
is given by evaluating elements of 
$A_m^*$ at $s_m$.
Since the sequence $\{s_m\,|\, m\geq 0\}$
is a sequence of divided powers, one may easily
check that
$\zeta$ is a character of $A^*$, 
i.e. $\zeta: A^*\to K$ is a homomorphism 
of K-algebras.
In terms of the notion introduced 
in \cite{ABS}, the pair $(A^*, \zeta)$ 
becomes {\it a combinatorial Hopf algebra}, 
and the homomorphism $\cS^*: A^* \to \cQf$
coincides the unique homomorphism guaranteed 
by Theorem $4.1$ in \cite{ABS} 
for combinatorial Hopf algebras.
\end{rmk}

\subsection{Possible Applications}\label{S2.3}

In this subsection, we discuss
the following possible applications 
of the universal property of 
the \cNcs system $(\cNsf, \Pi)$ 
given in Theorem \ref{Universal}. 
From the discussions below, 
we also can see some of 
the main motivations 
for the introduction of 
the \cNcs systems over 
associative algebras. 

First, the universal property of 
the \cNcs system $(\cNsf, \Pi)$ 
can  be used to solve any equations of 
Eqs.\,$(\ref{UE-1})$--$(\ref{UE-4})$
over any $K$-algebras $A$.
For example, given a $K$-algebra $A$ 
and $h(t)\in A[[t]]$ with $h(0)=0$, 
we can solve Eqs.\,$(\ref{UE-2})$ and 
$(\ref{UE-3})$ for $d(t)$ and 
$g(t)$ as follows. 
First, by Corollary \ref{Uniqueness-1}, 
we know, theoretically, 
there exists a unique \cNcs system 
$\Omega$ over $A$ with $h(t)$ as 
its fourth component. 
Hence,  by the universal 
property of $(\cNsf, \Pi)$,  
we have a unique homomorphism 
$\cS :\cNsf\to A$ such that 
$\cS^{\times 5} (\Pi)=\Omega$. 
On the other hand, 
the relations or polynomial identities 
between any two families 
of the NCSFs in the first four components of 
$\Pi$ have been given explicitly in \cite{G-T}. 
By applying the anti-involution 
$\omega_\Lambda$ in Proposition 
$\ref{omega-Lambda}$, 
one can easily derive the relations of 
the NCSFs $\Xi_m$'s with 
other NCSFs in $\Pi$ (for example, 
see $\S 4.1$ in \cite{GTS-III} 
for a complete list). 
In particular, the coefficients 
of $\sigma(t)$ 
and $\Phi(t)$ can be written as 
certain polynomials in 
the coefficients 
of $\psi(t)$.
Now, by simply applying the 
algebra homomorphism 
$\cS$ to these polynomials, 
we get the coefficients of the wanted 
solutions $d(t)$ and $g(t)$ 
in terms of the same polynomials 
in the coefficients of $h(t)$.
Hence, we get the solution 
$g(t)$ and $h(t)$ in $A[[t]]$. 
From the arguments above, 
we see that the generating 
functions of the NCSFs in the universal \cNcs system 
$(\cNsf, \Pi)$ can be viewed as 
the universal solutions to 
Eqs.\,$(\ref{UE-1})$--$(\ref{UE-4})$
over all associative $K$-algebras.

Secondly, suppose that 
a \cNcs system $(A, \Omega)$ 
is given. By applying the $K$-algebra 
homomorphism $\cS:\cNsf \to A$ 
guaranteed by the universal 
property of the system $(\cNsf, \Pi)$ 
to the identities of the NCSFs 
in the \cNcs system $\Pi$,
we get same identities for the corresponding 
elements of $A$ in the \cNcs system $\Omega$.  
This could be a very effective 
way to obtain identities for 
certain elements of $A$ if 
we could show that they are 
involved in a \cNcs system 
over $A$. For example, this gadget 
will be applied in the followed 
paper \cite{GTS-III} to derive some 
identities for certain differential 
operators which are important 
in the studies of the inversion 
problem (\cite{BCW}, \cite{E}), 
i.e. the problem to study 
various properties 
of the inverse maps of affine spaces.
On the other hand, 
if the given \cNcs system $(A, \Omega)$ 
has already been well-understood, 
the $K$-algebra homomorphism 
$\cS:\cNsf  \to A$ in turn 
gives a {\it specialization} 
or {\it realization} (\cite{G-T}, \cite{St2})
of NCSFs, which may be applied to study
certain properties of NCSFs.

\renewcommand{\theequation}{\thesection.\arabic{equation}}
\renewcommand{\therema}{\thesection.\arabic{rema}}
\setcounter{equation}{0}
\setcounter{rema}{0}

\section{\bf \cNcs Systems over Differential Operator Algebras}\label{S3}

In this section, we discuss 
the \cNcs systems that will be 
constructed in \cite{GTS-II} over 
differential operator algebras 
in commutative or noncommutative 
free variables. Certain properties 
of the resulting differential 
operator specializations 
of NCSFs, which will 
be proved in \cite{GTS-II} 
and \cite{GTS-V}, will also be discussed. 
The main purposes of this section and the 
next one are, first, to provide 
some supporting examples 
for the general discussions 
of \cNcs systems given 
in the previous section, 
and second, to give a shorter survey or 
review for some of the main results to be given 
in the followed papers \cite{GTS-II}, 
\cite{GTS-IV} and \cite{GTS-V}. For more examples 
of the specializations of NCSFs, see 
the references quoted in the introduction.

First, let us fix the following notation.

Let $K$ be any unital 
commutative $\bQ$-algebra as before
and $z=(z_1, z_2, ... , z_n)$ commutative 
or noncommutative free 
variables.\footnote{Since most of the results 
in this section do not depend 
on the commutativity 
of the free variables $z$,  
we will not distinguish the commutative 
and the noncommutative case, 
unless stated otherwise, 
and adapt the notations 
for noncommutative variables 
uniformly for the both cases.}
Let $t$ be a formal central parameter, 
i.e. it commutes with $z$ and elements of $K$.
We denote by $\kzz$ and 
$\kttzz$ the $K$-algebras of 
formal power series in $z$ over  
$K$ and $K[[t]]$, respectively.

By a {\it $K$-derivation} or simply {\it derivation} 
of $\kzz$, we mean a $K$-linear $\delta: \kzz\to \kzz$ 
that satisfies the Leibniz rule,
i.e. for any $f, g\in \kzz$, we have
 \begin{align}\label{Leibniz}
\delta (fg)=(\delta f)g+f(\delta g).
\end{align}
We will denote by 
$\cDrkzz$ or $\cDrzz$, when the base algebra
$K$ is clear from the context, 
the set of all $K$-derivations 
of $\kzz$.
The unital subalgebra of 
$\text{End}_k(\kzz)$ 
(endomorphisms of $\kzz$ as a $K$-vector space)
generated by all
$K$-derivations 
of $\kzz$ will be denoted by 
$\cDkzz$ or $\cDzz$. 
Elements of $\cDzz$ will be called 
{\it $(\text{formal})$ 
differential operators} in the
free variables $z$.

For any $\alpha\geq 1 $, 
we denote by $\cDrazz$ 
the set of the $K$-derivations 
of $\kzz$ which increase 
the degree in $z$ by 
at least $\alpha-1$. 
The unital subalgebra of 
$\cDzz$ generated by elements of 
$\cDrazz$ will be denoted by $\cDazz$. 
Note that, by the definitions above,
the operators of scalar multiplications
are also in $\cDzz$ and $\cDazz$.
When the base algebra is $K[[t]]$ 
instead of $K$ itself,
the notation 
$\cDrzz$, $\cDzz$, $\cDrazz$ 
and $\cDazz$ will be denoted by  
$\cDrtzz$, $\cDtzz$, $\cDrtazz$ 
and $\cDtazz$, respectively.
For example,  $\cDrtazz$ stands for 
the set of all $K[[t]]$-derivations 
of $\kttzz$, 
which increase the degree in $z$ 
by at least $\alpha-1$. 
Note that, $\cDrtazz=\cDrazz[[t]]$ and 
$\cDtazz=\cDazz[[t]]$.

For any $1\leq i\leq n$ and $u(z)\in \kzz$, 
we denote by $\lb u(z) \fr \p{\p z_i}\rb $
the $K$-derivation which maps $z_i$ to $u(z)$ and $z_j$ to $0$ 
for any $j\neq i$. 
For any $\Vec{u}=(u_1, u_2, \cdots, u_n)\in \kzz^{\times n}$, 
we set 
\begin{align}\label{Upz}
[\Vec{u}\pz]:=\sum_{i=1}^n [u_i \fr\p{\p z_i}].  
\end{align}

Note that, in the noncommutative case, 
we in general do {\bf not} have
$\lb u(z) \fr \p{\p z_i}\rb  g(z)  = u(z)  
\fr {\p g}{\p z_i}$ for all 
$u(z), g(z)\in \kzz$. This is the reason why
we put a bracket $[\cdot]$ in the notation above 
for the $K$-derivations.

With the notation above, 
it is easy to see that any 
$K$-derivations $\delta$ 
of $\kzz$ can be written uniquely  as 
$\sum_{i=1}^n \lb f_i(z)\fr\p{\p z_i}\rb$ 
with $f_i(z)=\delta z_i\in \kzz$ 
$(1\leq i\leq n)$.

With the commutator bracket, $\cDrazz$ $(\alpha\geq 1)$ 
forms a Lie algebra 
and its universal enveloping algebra is exactly 
the differential operator algebra
$\cDazz$.  Consequently, $\cDazz$ $(\alpha\geq 1)$
has a Hopf algebra structure as all 
other enveloping algebras of Lie algebras do.
In particular,
Its coproduct $\Delta$, antipode $S$ and co-unit $\epsilon$
are uniquely determined by the properties
\begin{align}
\Delta(\delta)&= 1\otimes\delta+\delta\otimes 1,\label{Coprd-delta} \\
S(\delta)&=-\delta, \label{antipd-delta}\\
\epsilon (\delta) &=\delta \cdot 1, \label{Counit-delta}
\end{align}
respectively, for any $\delta \in \cDrzz$.

For any $\alpha\geq 1 $, let $\ataz$ 
be the set of all the automorphisms 
$F_t(z)$ of $\kttzz$ over $K[[t]]$, 
which have the form
$F_t(z)=z-H_t(z)$ for some 
$H_t(z)\in \kttzz^{\times n}$ 
with $o(H_t(z))\geq \alpha$ 
and $H_{t=0}(z)=0$, where 
$o(H_t(z))$ denotes the minimum of 
the orders of all components of 
$H_t(z)$ as formal power series 
in $z$ with coefficients in $K[[t]]$.    
It is easy to check that $\ataz$ forms a 
subgroup of the automorphism group 
of $\kttzz$ over $K[[t]]$.
In particular, for any $F_t \in \ataz$ 
as above, its inverse map $G_t:=F_t^{-1}$
can always be written uniquely 
as $G_t(z)=z+M_t(z)$ \label{Mt(z)}
for some $M_t(z)\in \kttzz^{\times n}$ 
with $o(M_t(z))\geq \alpha$ 
and $M_{t=0}(z)=0$. 
Throughout this section, 
we will always let 
$H_t(z)$, $G_t(z)$ and $M_t(z)$
be determined as above. 

Now we fix an $\alpha \geq 1$ and 
an arbitrary $F_t\in \ataz$ and  
consider the \cNcs systems
$(\cDazz, \Oft)$ that will be 
constructed in \cite{GTS-II} over 
the differential operator algebra 
$\cDazz$.  Note that, $F_t \in \ataz$ 
can be viewed as a deformation 
parameterized by $t$ of 
the formal map $F(z)\!:=F_{t=1}(z)$, 
when it makes sense. 
For more studies on
$F_t\in \ataz$ from the deformation 
point view, see \cite{BurgersEq} 
and \cite{NC-IVP}. 
Actually, the construction of the 
\cNcs system $(\cDazz, \Oft)$ 
is mainly motivated by and 
also depends on the 
studies of $F_t\in \ataz$ given 
in \cite{BurgersEq} 
and \cite{NC-IVP}.

We first denote 
by the to-be-constructed
\cNcs system $\Omega_{F_t}$ as
\begin{align}\label{Def-Omega-Ft}
\Omega_{F_t}=(f(t),\, g(t),\, d(t),\, h(t),\, m(t)) 
\in \cDazz[[t]]^{\times 5}
\end{align}
and write the components of $\Oft$ above as in 
Eqs.\,(\ref{Def-lambda-m})--(\ref{Def-xi-m}) 
with the to-be-determined coefficients 
in $\cDazz$. Then the components of $\Oft$ 
are determined as follows.

The first three components of $\Oft$ are given 
by the following proposition which will be 
proved in Section $3.2$ in \cite{GTS-II}. 

\begin{propo}\label{TaylorExpansion-DLog}
There exist unique $f(t), g(t), d(t) \in \cDazz[[t]]$ 
with $f(0)=1$ and $d(0)=0$
such that, for any $u_t(z)\in \kttzz$, we have
\begin{align}
f(-t)\, u_t(z)&=u_t(F_t),\label{NewTaylorExpansion-e1}\\
g(t)\, u_t(z)&=u_t(G_t), \label{NewTaylorExpansion-e2}\\
e^{d(t)}\, u_t(z) &=u_t(G_t), \label{NewDLog-e2}
\end{align}
where, as usual, the exponential in 
Eq.\,$(\ref{NewDLog-e2})$ is given by 
\begin{align}\label{L2.2.3-e2}
e^{d(t)} = \sum_{m\geq 0} \frac {d(t)^m}{m!}.
\end{align}
\end{propo}

Note that, when we write $d(t)$ above as 
$d(t)= -\lb a_t(z)\pz \rb$ for some 
$a_t(z)\in t\kttzz$, then we get the so-called 
{\it D-Log}\label{D-Log} $a_t(z)$ of the automorphism 
$F_t(z)\in \ataz$, which has been studied 
in \cite{E1}--\cite{E3}, \cite{N}, \cite{Z-exp} 
and \cite{WZ} for the commutative case. 

The last two components of $\Oft$ are given directly as

\begin{align}
h(t)&:=\left[\frac{\p M_t}{\p t}(F_t) \pz \right ],\label{Def-h(t)}\\
m(t)&:=\left[\frac{\p H_t}{\p t}(G_t) \pz \right ]. \label{Def-m(t)}
\end{align}

To get some concrete ideas for the differential operators 
defined above, let us recall the following lemma proved 
in \cite{NC-IVP} for the special $F_t\in \ataz$ 
with $H_t(z)=tH(z)$ for some 
$H(z) \in \kzz^{\times n}$.

\begin{lemma} \label{NCIVP-L4.1.1}
For any $F_t \in \ataz$ of the form $F_t(z)=z-tH(z)$ 
as above, let $N_t(z)=t^{-1}M_t(z)$. 
Then we have

\begin{align}
m(t)&=\lb N_t(z)\pz \rb, \label{NCIVP-L4.1.1-e1} \\
h(t)&=\sum_{m\geq 1}  t^{m-1} \lb C_m(z)\pz \rb, \label{NCIVP-L4.1.1-e2}
\end{align}
where $C_m(z)\in \kzz^{\times n}$ $(m\geq 1)$ 
are defined recurrently by
\begin{align}
C_1(z)& =H(z),\label{NCIVP-L4.1.1-e3}\\
C_m(z)&=\lb C_{m-1}(z)\pz \rb H(z),\label{NCIVP-L4.1.1-e4}
\end{align}
for any $m\geq 2$.

Consequently, for any $m\geq 1$, the derivations 
$\psi_m$ and $\xi_m$ defined in 
Eqs.\,$(\ref{Def-m(t)})$ and $(\ref{Def-h(t)})$
are given by

\begin{align}
\psi_m &=\lb C_m(z)\pz \rb, \label{C4.1.2-e1}\\
\xi_m  &=\lb N_{[m]}(z)\pz \rb, \label{C4.1.2-e2}
\end{align}
where $N_{[m]}(z)\in \kzz^{\times n}$ $(m\geq 1)$ 
is the coefficient of $t^{m-1}$ of $N_t(z)$.
\end{lemma}

By the mathematical induction on $m\geq 1$, 
it is easy to show that, 
when $z$ are commutative variables, 
we further have
\begin{align}
C_m(z)=(JH)^{m-1}H(z) \label{Special-Cm}
\end{align}
for any $m\geq 1$, where $JH$ is the Jacobian matrix of 
$H(z)\in K[[z]]^{\times n}$.

\begin{theo} \label{S-Correspondence} $($\cite{GTS-II}$)$
For any $\alpha\geq 1 $ and $F_t(z)\in \ataz$, we have,

$(a)$ the $5$-tuple $\Omega_{F_t}$ defined as above
forms a \cNcs system over 
the differential operator algebra $\cDazz$.

$(b)$ let $(\cNsf, \Pi)$ be the \cNcs system
of NCSFs introduced in Section \ref{S2.2}, 
then there exists a unique homomorphism 
$\cS_{F_t}: \cNsf \to \cDazz$ of 
$K$-Hopf algebras such that 
$\cS_{F_t}^{\times 5}(\Pi)=\Omega_{F_t}$. 
\end{theo}

Note that, $(b)$ follows directly from $(a)$ and 
Theorem \ref{Universal}, since 
all the coefficients of $h(t)$ 
by Eq.\,(\ref{Def-h(t)}) 
are $K$-derivations and hence are 
primitive elements of 
the Hopf algebra $\cDrazz$.

For any $F_t(z)\in \ataz$, 
let $\Omega_{F_t}$ be defined above. 
We write the components of 
$\Omega_{F_t}$ as in Eq.\,(\ref{Def-Omega-Ft}) 
and coefficients of the components 
as in Eqs.\,(\ref{Def-lambda-m})--(\ref{Def-xi-m}). 
Then we have the following differential 
operator specializations of the NCSFs 
in the \cNcs system $(\cNsf, \Pi)$.

\begin{corol} \label{S-Correspondence-2}
For any $\alpha\geq 1 $ 
and $F_t(z)\in \ataz$, let 
$\cS_{F_t}\!:\! \cNsf\! \to \cDazz$ be 
the homomorphism of K-Hopf algebras 
in Theorem $\ref{S-Correspondence}$, $(b)$. 
Then, for any $m\geq 1$, we have 
 \begin{align}
 \mathcal S_{F_t} (\Lambda_m) & =\lambda_m, \label{L-l} \\
 \mathcal S_{F_t} (S_m) & =s_m, \label{S-s} \\
 \mathcal S_{F_t} (\Psi_m) & =\psi_m, \label{Psi-psi}  \\
 \mathcal S_{F_t} (\Phi_m)& =\phi_m, \label{Phi-phi}  \\
 \mathcal S_{F_t} (\Xi_m)& =\xi_m. \label{Xi-xi} 
 \end{align}
 \end{corol}

Note that, one direct consequence of 
Theorem \ref{S-Correspondence} above is 
the following well-defined map:
\begin{align}\label{Embedding}
\mathbb S: \ataz & \longrightarrow \text{\bf Hopf}\,(\cNsf, \cDazz)  \\
F_t \quad  & \longrightarrow \quad\quad\quad  \cS_{F_t}, \nno 
\end{align}
where $\text{\bf Hopf}\,(\cNsf, \cDazz)$ denotes 
the set of $K$-Hopf algebra 
homomorphisms from $\cNsf$ to $\cDazz$. 

Actually, as will be shown in \cite{GTS-II}, 
the following proposition holds. 

\begin{propo}\label{Isomorphism}
For any $\alpha \geq 1$, the map 
$\mathbb S$
defined in Eq.\,$(\ref{Embedding})$ is a bijection.
\end{propo}

Moreover, by identifying $\text{\bf Hopf}\,(\cNsf, \cDazz)$
with the set of all sequences of divided powers 
of the Hopf algebra
$\cDazz$, one can define a group product for the set
$\text{\bf Hopf}\,(\cNsf, \cDazz)$, with respect to 
which the bijection $\mathbb S$ above becomes 
an isomorphism of groups.

Next, let us consider the question when the  
Hopf algebra
homomorphism $\cS_{F_t}: \cNsf \to \cDazz$ 
preserves the gradings of $\cNsf$ and $\cDazz$.
Note that, precisely speaking, $\cDazz$ 
is not graded in the usual sense, 
for some infinite sums are allowed in $\cDazz$. 
But we can consider the following graded 
subalgebras of $\cDazz$. 

Let $\cDz$ be the differential operator 
algebra of the polynomial
algebra $\kz$, i.e. $\cDz$ is the unital subalgebra of 
$\text{End}_K(\kz)$ generated 
by all $K$-derivations of $\kz$.  
For any $m\geq 0$, let $\cD_{[m]}\langle z \rangle$ 
be the set of all differential operators 
$U$ such that, for any homogeneous polynomial 
$h(z)\in \kz$ of degree $d\geq 0$, $U h(z)$ 
either is zero or is homogeneous of degree $m+d$. 
For any $\alpha\geq 1$, set $\cDaz:=\cDz\cap \cDazz$.
Then, we have the grading
\begin{align}\label{Grading-cDz}
\cDaz &=\bigoplus_{m \geq \alpha-1} 
\cD_{[m]} \langle z \rangle, 
\end{align}
with respect to which $\cDaz$ 
becomes a graded $K$-Hopf algebra.

Now, for any $\alpha \geq 2$, we let 
$\gtaz$ be the set of all automorphisms 
$F_t\in \ataz$ such that $F_t(z)=t^{-1}F(tz)$
for some automorphism $F(z)$ of $\kzz$.
It is easy to check that $\gtaz$ 
is a subgroup of $\ataz$. 
Then we have the following proposition 
that will be proved in \cite{GTS-II}.

\begin{propo}\label{cS-graded}
For any $\alpha\geq 2$ and $F_t\in \ataz$, 
the differential operator specialization 
$\cS_{F_t}$ is a graded $K$-Hopf algebra 
homomorphism $\cS_{F_t}: \cNsf\to \cDaz\subset \cDazz$ 
iff $F_t\in \gtaz$.
\end{propo}

Now, for any $F_t\in \gtaz$ $(\alpha\geq 2)$, 
by the proposition above, 
we can take the graded dual of 
the graded $K$-Hopf algebra homomorphism 
$\cSft:\cNsf \to \cDaz$ 
and get the following corollary.

\begin{corol}\label{S-cQsf}
For any $\alpha \geq 2$ and $F_t\in\gtaz$, 
let $\cDaz^*$ be the graded dual 
of the graded $K$-Hopf algebra $\cDaz$. Then, 
$$
\cSft^*: \cDaz^* \to \cQf
$$ 
is a homomorphism of 
graded $K$-Hopf algebras. 
\end{corol}

Next, let us point out 
the following property to be 
proved in \cite{GTS-V}
of the differential operator 
specializations 
$\cSft$ $(F_t\in \ataz)$.

For any  $\alpha\geq 1$,
let $\mathbb B^{[\alpha]}_t\langle z \rangle$
be the set of automorphisms $F_t=z-H_t(z)$ 
of the polynomial algebra 
$\ktz$ over $K[t]$ such that 
the following conditions are satisfied.
\begin{enumerate}
\item[$\bullet$]  $H_{t=0}(z)=0$. 
\item[$\bullet$] 
$H_t(z)$ is homogeneous in $z$ 
of degree $d \geq \alpha$.
\item[$\bullet$] 
With a proper permutation of the free variables $z_i$'s,
the Jacobian matrix $JH_t(z)$ becomes strictly lower triangular.
\end{enumerate}

\begin{theo}\label{StabInjc-best} 
In both commutative and noncommutative cases, 
the following statement holds. 
\vskip2mm
For any fixed $\alpha\geq 1$ and non-zero NCSF 
$P \in \cNsf$, there exist $n\geq 1$ 
$($the number  of the free variable 
$z_i$'s$)$ and $F_t(z)\in \mathbb B^{[\alpha]}_t\langle z\rangle$ 
such that $\cS_{F_t} (P)\neq 0$.
\end{theo}

\begin{rmk}
As pointed out earlier 
in Subsection \ref{S2.3}, 
we can apply the homomorphism 
$\cSft:\cNsf\to \cDazz$ 
to transform the identities 
of the NCSFs to 
the identities of the corresponding 
differential operators 
in the \cNcs systems $\Oft$. 
Combining with the special forms 
of the differential operators 
$\psi_m$'s and $\xi_m$'s in 
Eqs.\,$(\ref{C4.1.2-e1})$ and 
$(\ref{C4.1.2-e2})$, 
respectively,
we can derive more 
identities for 
the inverse map, the D-Log of 
$F_t\in \ataz$ as well as 
the formal flow generated 
by $F_t(z)$, which may be  
applied further to study the inversion problem. 
For detailed discussions in this 
direction, see the sequel paper 
\cite{GTS-III}. 
\end{rmk}

Finally, let us summarize 
the main results discussed 
in this section as follows.
By Theorem \ref{S-Correspondence}, 
for any $F_t \in\ataz$, we have a specialization 
$\cSft: \cNsf \to \cDazz$ of NCSFs by 
differential operators in $\cDazz$; 
by Proposition \ref{Isomorphism}, we know
any such a specialization of NCSFs, if it 
is also a $K$-Hopf algebra homomorphism,  
is give by $\cSft$ for some $F_t \in \ataz$; 
By Proposition \ref{cS-graded}, we know
exactly when the specialization 
$\cSft:\cNsf\to \cDaz$ 
preserves the gradings of $\cNsf$ and $\cDaz$;
Finally, by Theorem \ref{StabInjc-best},  we know
the smaller family of the specializations $\cSft$ 
with all possible $n\geq 1$ and 
$F_t\in \mathbb B^{[\alpha]}_t\langle z\rangle$
is already fine enough to distinguish any two 
different NCSFs.

\renewcommand{\theequation}{\thesection.\arabic{equation}}
\renewcommand{\therema}{\thesection.\arabic{rema}}
\setcounter{equation}{0}
\setcounter{rema}{0}

\section{\bf  A \cNcs System over 
the Grossman-Larson Hopf Algebra 
of Labeled Rooted Trees} 
\label{S4}

In this section, we fix a non-empty 
$W\subseteq \bN^+$ and
first recall the Connes-Kreimer Hopf algebra 
$\cH_{CK}^W$ (\cite{CM}, \cite{Kr}, \cite{CK}, \cite{F}) 
and the Grossman-Larson Hopf algebra 
$\mathcal H^W_{GL}$ (\cite{GL}, \cite{CK}, \cite{F})
of $W$-labeled rooted forests and 
$W$-labeled rooted trees, respectively.
We then discuss the 
\cNcs system $(\cH_{GL}^W, \Omega_\bT^W)$ 
that will be constructed in \cite{GTS-V} over 
the Grossman-Larson Hopf algebra 
$\mathcal H^W_{GL}$ and certain properties 
of the resulting specializations of NCSFs 
by $W$-labeled rooted trees. 
Finally, we briefly 
explain a connection,
which will be given in \cite{GTS-V}, 
between the \cNcs system 
$(\cH_{GL}^W, \Omega_\bT^W)$ 
with the \cNcs system $(\cDazz, \Oft)$ 
$(F_t\in \ataz)$ discussed in Section \ref{S3}.
Some consequences of this connection 
will also be discussed.

First, let us fix the following notation.

\vskip2mm

{\bf Notation:}

\vskip2mm

By a {\it rooted tree} we mean a finite
1-connected graph with one vertex designated as its {\it root.} 
For convenience, we also view the empty set $\emptyset$ 
as a rooted tree and call it the {\it emptyset} rooted tree.
The rooted tree with a single vertex 
is called the {\it singleton} 
and denoted by $\circ$. 
There are natural ancestral relations between vertices.  
We say a vertex $w$ is a {\it child} of vertex $v$ 
if the two are connected by an
edge and $w$ lies further from the root than $v$.  
In the same situation, we say $v$ is the {\it parent} of $w$.  
A vertex is called a {\it leaf}\/ if it has no
children. 

Let $W\subseteq \bN^+$ be 
any non-empty subset of positive integers. 
A {\it $W$-labeled rooted tree}
is a rooted tree 
with each vertex labeled by 
an element of $W$. If an element $m \in W$ 
is assigned to a vertex $v$, 
then $m$ is called the {\it weight} 
of the vertex $v$. When we speak of isomorphisms 
between unlabeled (resp.\,\,$W$-labeled) rooted trees, 
we will always mean isomorphisms 
which also preserve 
the root (resp.\,\,the root and also the labels of vertices).
We will denote by $\mathbb T$ (resp.\,\,$\mathbb T^W$) 
the set of isomorphism classes 
of all unlabeled (resp.\,\,$W$-labeled) rooted trees. 
A disjoint union of any finitely many rooted trees 
(resp.\,\,$W$-labeled rooted trees) 
is called a {\it rooted forest} (resp.\,\,$W$-labeled {\it rooted forest}).
We denote by $\mathbb F$ (resp.\,\,$\mathbb F^W$) 
the set of unlabeled (resp.\,\,$W$-labeled) rooted forests.  

\vskip3mm

With these notions in mind, we establish the following notation.
\begin{enumerate}
\item For any rooted tree $T\in \bT^W$, we set the following notation:
\begin{itemize}
\item $\text{rt}_T$ denotes the root vertex of $T$ and $O(T)$
the set of all the children of $\text{rt}_T$. We set
$o(T)=|O(T)|$ (the cardinal number of the set $O(T)$). 

\item $E(T)$ denotes the set of edges of $T$.

\item $V(T)$ denotes the set of  vertices of $T$ and $v(T)=|V(T)|$.

\item $L(T)$ denotes the set of leaves of $T$ and $l(T)=|L(T)|$


\item For any $T\in \bT^W$ and $T\neq \emptyset$, 
$|T|$ denotes the sum of the weights of 
all vertices of $T$. 
When $T=\emptyset$, we set $|T|=0$.

\item For any $T\in \bT^W$, we denote 
by $\text{Aut}(T)$ the automorphism group 
of $T$ and $\alpha(T)$ the cardinal 
number of $\text{Aut}(T)$.

\item For any $v\in V(T)$, we define the {\it height} of $v$ 
to be the number of edges in the (unique) geodesic connecting 
$v$ to $\text{rt}_{T}$.
The {\it height} of $T$ is defined to be the maximum of the heights of
its vertices.

 \end{itemize}

\item Any subset of $E(T)$ is called a {\it cut} of $T$. 
A cut $C\subseteq E(T)$ is said to be {\it admissible} 
if no two different edges of $C$ lie in 
the path connecting the root and a leaf. We denote by 
$\mathcal C(T)$ the set of all admissible cuts of $T$. 
Note that, the empty subset $\emptyset$ of $E(T)$
and $C=\{e\}$ for any $e\in E(T)$ 
are always admissible cuts. 

\item For any 
$T \in \bT^W$ with $T\neq \circ$, 
let $C\in \mathcal C(T)$ 
be an admissible cut of $T$ 
with $|C|=m\geq 1$.
Note that, after deleting 
the edges in $C$ from $T$, 
we get a disjoint union of $m+1$ 
rooted trees, 
say $T_0$, $T_1$, ..., $T_m$ 
with $\text{rt}(T)\in V(T_0)$.
We define $R_C(T)=T_0 \in \mathbb T^W$ 
and $P_C (T)\in \mathbb F^W$ 
the rooted forest formed by $T_1$, ..., $T_m$. 

\item For any $T\in \bT^W$, we say $T$ 
is a {\it chain} if 
its underlying rooted tree
is a rooted tree with a single leaf.
We say $T$ is a {\it shrub} if 
its underlying rooted tree
is a rooted tree 
of height $1$. 
We say $T$ is {\it primitive} if 
its root has only one child. 
For any $m\geq 1$, we set $\mathbb H_m$, 
$\mathbb S_m$ and $\mathbb P_m$ to be the sets of 
the chains, shrubs and primitive rooted trees 
$T$ of weight $|T|=m$, respectively. 
$\mathbb H$, $\mathbb S$ and $\mathbb P$ 
are set to be the unions of $\mathbb H_m$, 
$\mathbb S_m$ and $\mathbb P_m$, 
respectively, for all $m\geq 1$.   

\end{enumerate}

Let $K$ be any unital commutative $\bQ$-algebra 
and $W$ a non-empty subset of positive integers.
First, let us recall the Connes-Kreimer Hopf algebras 
$\mathcal H_{CK}^W$ of labeled rooted forests. 

As a $K$-algebra, the Connes-Kreimer Hopf algebra 
$\mathcal H_{CK}^W$ is the free commutative algebra 
generated by formal variables 
$\{X_T \,|\, T\in \mathbb T^W \}$. 
Here, for convenience, we will still use $T$ 
to denote the variable $X_T$ in 
$\mathcal H_{CK}^W$.
The $K$-algebra product is given by the disjoint union. 
The identity element of this algebra, denoted by $1$,  
is the free variable $X_\emptyset$ 
corresponding to the emptyset rooted tree. 
The coproduct $\Delta: \mathcal H_{CK}^W \to 
\mathcal H_{CK}^W \otimes \mathcal H_{CK}^W$ 
is uniquely determined by setting 
\begin{align}
\Delta(1)&=1\otimes 1, \label{CK-Delta-1} \\
\Delta(T)&= T\otimes 1+ 
\sum_{C\in \mathcal C(T)} P_C(T) \otimes R_C(T). \label{CK-Delta-2}
\end{align}

The co-unit $\epsilon: \mathcal H_{CK}^W \to K$ 
is the $K$-algebra homomorphism 
which sends $1\in \mathcal H_{CK}^W$ to $1\in K$ 
and $T$ to $0$ for any $T\in \mathbb T^W$ 
with $T\neq \emptyset$. 
With the operations defined above and 
the grading given by the weight, 
the vector space $\mathcal H_{CK}^W$ 
forms a connected graded commutative 
bi-algebra. Since any connected graded bialgebra 
is a Hopf algebra, there is a unique antipode 
$S: \mathcal H_{CK}^W\to \mathcal H_{CK}^W$ 
that makes $\mathcal H_{CK}^W$ a connected 
graded commutative $K$-Hopf algebra. 
For a formula for the antipode, see \cite{F}.

Next we recall the Grossman-Larson Hopf algebra 
of labeled rooted trees. 
First we need define the following operations 
for labeled rooted forests.  For any 
labeled rooted forest $F$ which is disjoint 
union of labeled 
rooted trees $T_1$, $T_2$, ... , $T_m$, 
we set $B_+(T_1, T_2, \cdots, T_m)$ 
the rooted tree obtained by connecting roots 
of $T_i$ $(1\leq i\leq m)$ 
to a newly added root. We will keep the labels 
for the vertices of $B_+(T_1, T_2, \cdots, T_m)$ 
from $T_i$'s, but for the root, we label it by $0$. 
 
Now, we set $\bar {\mathbb T}^W:=\{ B_+(F) \, | \, F \in \mathbb F^W \}$.
Then,  $B_+: \mathbb F^W \to \bar{\mathbb T}^W$ becomes a bijection.
We denote by $B_- :  \bar{\mathbb T}^W \to \mathbb F^W$ 
the inverse map of $B_+$. More precisely, for any 
$T\in \bar{\mathbb T}^W$, 
$B_-(T)$ is the $W$-labeled rooted forest obtained by cutting off 
the root of $T$ as well as all edges connecting to the root in $T$.

Note that, precisely speaking,
elements of $\bar{\mathbb T}^W$ are not 
$W$-labeled trees for $0\not \in W$. 
But, if we set 
$\bar W=W\cup\{0\}$, 
then we can view $\bar{\mathbb T}^W$ 
as a subset 
of $\bar W$-labeled 
rooted trees $T$ with the root $\text{rt}_T$ 
labeled by $0$ 
and all other vertices 
labeled by non-zero elements of $\bar W$. 
We extend the definition of
the weight for elements 
of $\mathbb F^W$ to elements of 
$\bar{\mathbb T}^W$ by simply counting 
the weight of roots by zero. 
We set $\bar{\mathbb S}_m^W:=B_+(\mathbb S_m^W)$
$(m\geq 1)$ and $\bar{\mathbb S}^W:=B_+(\mathbb S^W)$. 
We also define
$\bar{\mathbb H}_m^W$, 
$\bar{\mathbb P}_m^W$, 
$\bar{\mathbb H}^W$ and 
$\bar{\mathbb P}^W$ 
in the similar way.

The Grossman-Larson Hopf algebra $\mathcal H_{GL}^W$ as a vector space 
is the vector space spanned by elements of $\bar{\mathbb T}^W$ 
over $K$. For any $T\in \bar{\mathbb T}^W$, we will still denote by 
$T$ the vector in $\mathcal H_{GL}^W$ that is corresponding to $T$. 
The algebra product is defined as follows. 
For any $T, S\in \bar{\mathbb T}^W$ with 
$T=B_+(T_1, T_2, \cdots, T_m)$, we set $T\cdot S$ to be the sum of  
the rooted trees obtained by connecting the roots of $T_i$ 
$(1\leq i\leq m)$ to vertices of $S$ 
in all possible $m^{v(S)}$ different ways. 
Note that, the identity element with respect to this 
algebra product is given by the singleton 
$\circ=B_+(\emptyset)$. But we will denote it by $1$.

To define the co-product $\Delta: \cH_{GL}^W \to \cH_{GL}^W \otimes \cH_{GL}^W$,
we first set 
\begin{align}
\Delta (\circ)=\circ \otimes \circ. 
\end{align}
 
Now let $T\in \bar{\mathbb T}^W$ with $T\neq \circ$, 
say $T=B_+ (T_1, T_2, \cdots, T_m)$ with
$m \geq 1$ and $T_i\in \mathbb T^W$ 
$(1\leq i\leq m)$. 
For any non-empty subset 
$I\subseteq \{1, 2, \cdots, m\}$,
we denote by 
$B_+(T_I)$ the rooted tree 
obtained by applying the $B_+$ operation 
to the rooted trees $T_i$ 
with $i\in I$. 
For convenience, when $I=\emptyset$, 
we set $B_+(T_I)=1$. 
With this notation fixed, the co-product 
for $T$ is given by
\begin{align}
\Delta (T)=\sum_{I\sqcup J=\{1, 2, \cdots, m\}} 
B_+(T_I) \otimes B_+(T_J).
\end{align}

Note that, a rooted tree in 
$\bar{\mathbb T}^W$ is a primitive element of 
the Hopf algebra $\cH_{GL}^W$ 
iff it is a primitive rooted tree 
in the sense that we defined before, 
namely the root of $T$ has 
one and only one child.

The co-unit $\epsilon: \mathcal H_{GL}^W \to K$ 
is the $K$-algebra homomorphism 
which sends $1\in \mathcal H_{CK}^W$ to $1\in K$ 
and $T$ to $0$ for any $T\in \mathbb T^W$ 
with $T\neq \emptyset$. With the operations defined 
above and the grading given by the weight, 
the vector space $\mathcal H_{GL}^W$ 
forms a connected graded commutative 
bi-algebra. Since any connected graded bialgebra 
is a Hopf algebra, there is a unique antipode 
$S: \mathcal H_{GL}^W\to \mathcal H_{CK}^W$ 
that makes $\mathcal H_{GL}^W$ a connected 
graded commutative $K$-Hopf algebra. 
For a formula for the antipode, 
see \cite{GTS-IV}.

Note that it has been shown in \cite{H} and \cite{F} that,
the Grossman-Larson Hopf algebra $\cH_{GL}^W$ and 
the Connes-Kreimer $\cH_{CK}^W$ are graded dual 
to each other.
The pairing is given by, 
for any $T\in \bar{\mathbb T}^W$ 
and $S\in \mathbb F^W$, 
\begin{align}\label{pairing}
<T, F>=\begin{cases} 0, & \text{ if } T \not \simeq B_+(F),\\
\alpha (T), &\text{ if } T \simeq B_+(F).
\end{cases}
\end{align}

Now we consider the \cNcs system $\Omega_{\bT}^W$ that will be 
constructed in \cite{GTS-IV} over the Grossman-Larson Hopf algebra 
$\cH_{GL}^W$.

First, let us define the following constants for the rooted trees in $\bbT^W$:
\begin{enumerate}
\item[$\bullet$] We set $\beta_T$ to be the weight 
of the unique leaf of $T$ if $T\in \bar{\mathbb H}^W$ and $0$ otherwise. 
\item[$\bullet$] We set $\gamma_T$ to be the weight 
of the unique child of the root of $T$ if $T\in \bar{\mathbb P}^W$ and $0$ otherwise. 
\item[$\bullet$]\label{theta} We set $\theta_T$ to be the coefficient of $s$ 
of the order polynomial $\Omega(B_-(T), s)$ of the underlying 
unlabeled rooted forest of $B_-(T)$.
\end{enumerate}

For general studies on the order 
polynomials $\Omega(P, s)$ 
of finite posets $P$, see \cite{St1}. 
For a combinatorial interpretation of the constant
$\phi_T:=(-1)^{v(T)-1}\theta_{B_+(T)}$  
in terms of the numbers of chains 
with fixed lengths in  
the lattice of the ideals of the poset $T$, 
see Lemma $2.8$ in \cite{SWZ}.

Now we consider the following generating functions of $T\in \bbT^W$.

\allowdisplaybreaks{
\begin{align}
\tilde f(t):&=\sum_{T \in \bar{\mathbb S}^W} (-1)^{o(T)-|T|}  t^{|T|} {\mathcal  V_T}
=1 + \sum_{\substack{T \in \bar{\mathbb S}^W \\
T \neq \circ }} (-1)^{o(T)-|T|}  t^{|T|} {\mathcal  V_T}, 
\label{Def-3f(t)} \\ 
\tilde g(t):&= \sum_{T \in \bar{\mathbb T}^W} t^{|T|} {\mathcal  V_T}
=1 + \sum_{\substack{ T  \in \bar{\mathbb T}^W \\
T \neq \circ }} t^{|T|} {\mathcal  V_T}, 
\label{Def-3g(t)}\\ 
\tilde d(t):&= \sum_{T \in \bar{\mathbb P}^W}  t^{|T|} \theta_T {\mathcal  V_T}. 
\label{Def-3d(t)}\\
\tilde h(t):&= \sum_{T \in \bar{\mathbb H}^W} t^{|T|-1}  \beta_T \mathcal  V_T, 
\label{Def-3h(t)} \\ 
\wtilde m(t):&= \sum_{T \in \bar{\mathbb P}^W}  t^{|T|-1} \gamma_T {\mathcal  V_T}, 
\label{Def-3m(t)}
\end{align} 
where, for any $T\in \bbT^W$, $\mathcal V_T:=\frac 1{\alpha(T)} T.$
}
We further set
\begin{align}\label{Def-3-Omega}
\Omega_\bT^W:=
(\, \tilde f(t),\, \tilde g(t), 
\, \tilde d\,(t), \, \tilde h(t), \wtilde m(t)\, ). 
\end{align}

\begin{theo}\label{Main-Thm-Trees} $($\cite{GTS-IV}$)$\,
 For any non-empty set $W\subseteq\bN^+$, we have

$(a)$ the $5$-tuple $\Omega_\bT^W$ defined 
in Eq.\,$(\ref{Def-3-Omega})$
forms a \cNcs system over 
the Grossman-Larson Hopf algebra $\mathcal H^W_{GL}$.

$(b)$ let $(\cNsf, \Pi)$ be the \cNcs system
of NCSFs introduced in Section \ref{S2.2}, 
then there exists a unique graded 
$K$-Hopf algebra homomorphism 
$\cT_W: \cNsf \to \mathcal H^W_{GL}$ such that 
$\cT^{\times 5}_W(\Pi)=\Omega_\bT^W$. 
\end{theo}

Note that the grade duals of $\cNsf$ and $\cH_{GL}^W$ are 
the graded $K$-Hopf algebras $\cQf$ of 
quasi-symmetric functions and the Connes-Kreimer 
Hopf algebra $\cH_{CK}^W$, respectively.
Since the $K$-Hopf algebra homomorphism 
$\cT_W: \cNsf\to \cH_{GL}^W$ in Theorem \ref{Main-Thm-Trees}
preserves the gradings, by taking the graded duals, we have  
the following correspondence. 

\begin{corol}\label{cQsf}
For any non-empty $W\subseteq \bN^+$, $\mathcal T^*_W: \cH_{CK}^W\to \cQf$ is 
a homomorphism of graded $K$-Hopf algebras.  
\end{corol}

Furthermore, the following result will be proved in \cite{GTS-V}.

\begin{propo}\label{Injective}
When $W=\bN^+$, the graded Hopf algebra homomorphism 
$\cT_W: \cNsf \to \cH_{GL}^W$ in Theorem \ref{Main-Thm-Trees} 
is an embedding. Consequently, its graded dual 
$\mathcal T^*_W: \cH_{CK}^W \to \cQf$ in this case is a 
surjective homomorphism of graded Hopf algebras.
\end{propo}

\begin{rmk}\label{R4.8}
As we mentioned earlier in Subsection \ref{S2.3},
by applying the specialization 
${\mathcal T}_W:\cNsf \to \cH_{GL}^W$ in Theorem \ref{Main-Thm-Trees},  
we will get a host of identities from 
the identities of the NCSFs in the \cNcs system 
$(\cNsf, \Pi)$ for the $W$-rooted trees 
in the \cNcs system $(\cH_{GL}^W, \Omega_\bT^W)$.
We believe some of these identities 
are interesting from the aspect of combinatorics 
of rooted trees. 
But, to keep this paper in a certain size, 
we have to ask the reader who is interested in 
these identities to do the translations
via the Hopf algebra homomorphism 
${\mathcal T}_W:\cNsf \to \cH^W_{GL}$.
\end{rmk}

Next, we briefly explain a connection,
which will be given in \cite{GTS-V}, 
between the \cNcs system $\Omega_\bT^W$ over the 
Grossman-Larson Hopf algebra $\cH_{GL}^W$ 
and the \cNcs system $\Oft$  
$(F_t\in \ataz)$ 
over the differential operator algebra $\cDazz$
discussed in Section \ref{S3}.
This connection will play an important role 
in the proofs of Theorem \ref{StabInjc-best} 
and Proposition \ref{Injective} mentioned before.

Let $W\subseteq \bN^+$ be any non-empty subset of 
positive integers and 
$F_t=z-H_t(z) \in \ataz$ such that $H_t(z)$ 
can be written as $\sum_{m\in W} t^m H_{[m]}(z)$ 
for some $H_{[m]}(z)\in \kzz^{\times n}$ $(m\in W)$.
In \cite{GTS-V},
a Hopf algebra homomorphism 
$\cAft: \cH_{GL}^W \to \cDazz$
such that $\cAft^{\times 5}(\OTW)=\Oft$ 
will be constructed. 
Furthermore, with this Hopf algebra homomorphism 
$\cAft$, we have the following 
commutative diagrams. 

\begin{propo}\label{CD-2}
For any $\alpha \geq 1$, let $W \subseteq\bN^+$ and 
$F_t\in \ataz$ fixed as above, 
we have the following commutative diagrams of 
$K$-Hopf algebra homomorphisms.
\begin{align}
\begin{CD}\label{CD-2-e1}
\cNsf @>{\cT_W} >> {\mathcal H}_{GL}^W \\
@V \cS_{F_t} VV @V \mathcal A_{F_t}  VV\\
\cDazz @= \cDazz
\end{CD}
\end{align}
\end{propo}

Combining Proposition \ref{S-cQsf} and 
Proposition \ref{CD-2} above, 
we have the following proposition.

\begin{propo}\label{CD-3}
For any $\alpha \geq 2$ and $F_t\in\gtaz$, 
we have the following commutative diagrams of $K$-Hopf algebra homomorphisms.
\begin{align}
\begin{CD}
{\mathcal Q}Sym @ < {{\mathcal T}_W^*} << {\mathcal H_{CK}^W } \\
@A {\cS_{F_t}^*} AA @A {\mathcal A^*_{F_t} } AA  \\
{\cDaz^* } @= {\cDaz^*}
\end{CD}
\end{align}
\end{propo}

Finally, let us point out that,
by applying the Hopf algebra homomorphism 
$\mathcal A_{F_t}: \cH_{GL}^W \to \cDazz$ 
above, we can also derive 
the tree expansion formulas 
for the inverse map $G_t(z)$, the D-Log $a_t(z)$ 
and the formal flow generated by $F_t\in \ataz$, 
which will generalize the tree expansion formulas 
obtained in \cite{BCW}, \cite{Wr3} and \cite{WZ} in 
the commutative case to the noncommutative case. 
For more details, see \cite{GTS-V}.

{\small \sc Department of Mathematics, Illinois State University,
Normal, IL 61790-4520.}

{\em E-mail}: wzhao@ilstu.edu.

\end{document}